\documentclass[]{informs3}

\OneAndAHalfSpacedXI

\usepackage{natbib}
 \bibpunct[, ]{(}{)}{,}{a}{}{,}%
 \def\bibfont{\small}%

\usepackage{fix-cm}

\usepackage{tikz}
\usepackage{pgfplots}
\usepackage{graphicx}
\usepackage{subcaption}
\graphicspath{ {./images/} }
\usepackage{amsmath}
\usepackage{multirow}
\usepackage{diagbox}
\usepackage{float}

\usepackage{multirow}
\usepackage{algorithm}
\usepackage[noend]{algorithmic}

\usepackage{array}

\usepackage{tikz,pgfplots}
\usetikzlibrary{matrix}
\usepgfplotslibrary{groupplots}
\pgfplotsset{compat=newest}
\usetikzlibrary{calc, positioning, shapes, fit, arrows, shadows, chains}
\tikzstyle{line} = [ draw, -latex']
\usetikzlibrary{shapes, arrows}
\usepgfplotslibrary{units}
\usepackage{url}

\usepackage{booktabs}

\usepackage{color}

\definecolor{color1}{rgb}{1.0, 0.5, 0.0}
\definecolor{color2}{rgb}{0.9, 0.4, 0.3}
\definecolor{color3}{rgb}{0.8, 0.3, 0.5}
\definecolor{color4}{rgb}{0.6, 0.2, 0.8}
\definecolor{color5}{rgb}{0.4, 0.1, 0.9}
\definecolor{color6}{rgb}{0.2, 0.0, 1.0}

\begin{document}

\RUNAUTHOR{Kayac\i k et al.}
\RUNTITLE{Dual Sourcing of Green Hydrogen}

\TITLE{Dual Sourcing of Green Hydrogen: Balancing Local Production with Stochastic Capacity and Import with Random Yield}

\ARTICLEAUTHORS{%
\AUTHOR{Sezen Ece Kayac{\i}k}
	\AFF{Department of Operations, Faculty of Economics and Business, University of Groningen, the Netherlands, \EMAIL{s.e.kayacik@rug.nl}}
\AUTHOR{Albert H. Schrotenboer}
	\AFF{Operations, Planning, Accounting and Control Group, School of Industrial Engineering, Eindhoven University of Technology, the Netherlands,  \EMAIL{a.h.schrotenboer@tue.nl}}
\AUTHOR{Iris F. A. Vis}
\AFF{Department of Operations, Faculty of Economics and Business, University of Groningen, the Netherlands, \EMAIL{i.f.a.vis@rug.nl}}
\AUTHOR{Beste Basciftci}
\AFF{Department of Business Analytics, Tippie College of Business, University of Iowa, Iowa City, Iowa, the United States, \EMAIL{beste-basciftci@uiowa.edu}}
\AUTHOR{Evrim Ursavas}
\AFF{Department of Operations, Faculty of Economics and Business, University of Groningen, the Netherlands, e.ursavas@rug.nl}
 }
 
\ABSTRACT{
Green hydrogen is a critical component for achieving the European Union’s 2050 net-zero emissions goal. However, ensuring a reliable and stable supply is challenging, particularly when local production of green hydrogen is subject to high variability due to fluctuating renewable energy output. Although importing from regions with stable renewable resources may offer greater reliability, they introduce longer lead times and potential energy losses during transportation and conversion. To address this issue, we develop optimal dual sourcing policies for green hydrogen through modeling and solving a Markov Decision Process that integrates general lead times, stochastic local supply capacity, and random yield from import. Alongside optimal dual sourcing policies, we propose heuristic policies that offer both flexibility and stability, enabling practical implementation while achieving near-optimal performance.  We test our framework on a case study for the Netherlands and obtain the following insights: (i) based on the results across different countries and cost settings, our dual sourcing model demonstrates an average cost benefit of 8\% compared to models that ignore stochastic supply capacity and random yield; 
(ii) the proposed heuristic policies can perform comparably to optimal policies under varying country-specific conditions and cost settings, offering insights for shaping hydrogen trade agreements between importing and exporting countries; (iii) sensitivity analyses on production and storage costs, as well as variability in supply capacity, demand, and random yield, reveal the conditions needed to achieve specific local production levels. These results thereby support feasibility of the Netherlands' climate scenarios and ambitions while guiding green hydrogen investment planning.
}


\maketitle

\vspace{-0.5cm}

\section{Introduction}

The European Union aims to achieve net zero emission targets by 2050 in accordance with global climate mitigation goals \citep{eu2012}. In this context, green hydrogen has emerged as a promising energy source, offering a versatile and sustainable solution to mitigate carbon emissions in various sectors such as mobility, industry, and heating \citep{staffell2019role}. Green hydrogen is clean energy that is produced through water electrolysis, using electricity from renewable sources like solar or wind \citep{hosseini2016hydrogen}. 
However, in countries with high variability in renewable energy production, relying solely on local production poses a challenge in maintaining a reliable green hydrogen supply. 
A solution is to import hydrogen from countries with substantial renewable energy potential, such as Morocco and Norway, given the transportability of green hydrogen and its derivatives such as green ammonia \citep{IRENA2021GreenHydrogen}. 
However, importing hydrogen introduces additional challenges, such as longer lead times and energy loss during transportation and conversion processes \citep{IRENA2022Geopolitics}. Additionally, local production and import cost dynamics can vary significantly across countries. Given these complexities, local production and import must be optimally coordinated to ensure a reliable and cost-efficient green hydrogen supply.



International collaboration is expected to play a significant role in matching supply and demand for green hydrogen \citep{McKinsey2023Hydrogen}. Currently, numerous governments are already establishing bilateral agreements to facilitate international hydrogen trade. Examples include partnerships between the Netherlands and Morocco, Germany and Chile, South Korea and Saudi Arabia, as well as Japan and Australia \citep{IRENA2022Geopolitics}. 
Looking ahead, in 2050, approximately one-quarter of the world's total hydrogen demand is expected to be met through international trade \citep{IRENA2022Hydrogen}. As these international partnerships grow, the specifics of trade agreements—such as order quantities and delivery schedules—will need to be carefully negotiated and structured. A key challenge is that exporting countries need stability in production quantities and schedules for effective planning, while importing countries face local supply uncertainties and fluctuating demand, making it difficult to commit to fixed quantities. Consequently, balancing stability and flexibility is critical to ensure that supply mechanisms can adapt to changing conditions while maintaining practicability \citep{basciftci2019adaptive, kayacik2024partially}. 

This paper is the first to study optimal dual sourcing policies between local production and import of green hydrogen. We consider a dual sourcing inventory control model, where local production and import decisions are dynamically adjusted. 
Our objective is to find a cost-minimizing policy that dynamically prescribes, based on the system's current state, the optimal quantities to source from each supply option. We consider that local production operates with stochastic capacity, meaning that its output fluctuates randomly over time. In addition, we consider random yield for importing, meaning that only a variable proportion of the actual order is delivered. This refers to the uncertainty in the amount of hydrogen that arrives, for example, due to potential losses during transportation and conversion processes. We assume that the import lead time is strictly greater than that of local production. 

We model the system as a Markov Decision Process (MDP) considering general lead times, stochastic capacity of the local supplier, and random yield of import. Via the value iteration method, we obtain optimal dual sourcing policies. Inspired by practice, we develop four heuristic dual sourcing policies that offer different levels of flexibility to handle variability in order quantities while ensuring stability for practical implementation. These policies support practitioners in forming effective hydrogen trade agreements under complex supply chain dynamics.

We apply our framework to a case study in the Netherlands based on real-world data. In our case study, we first show the cost benefits of the dual sourcing model with stochastic supply capacity and random yield, comparing it to single-sourcing models and dual sourcing models that ignore stochastic supply capacity, random yield, and both. We show that considering both stochastic supply capacity and random yield demonstrates an average cost benefit of 8\%. We then test the proposed heuristic dual sourcing policies to analyze the conditions in which they perform close to optimal policies in terms of cost efficiency.  
Our experiments show that increasing flexibility in ordering, from fixed amounts to partially adjustable policies, improves cost efficiency by an average of 11\%, and performs within 2\% of the optimal policy. In turn, our analysis provides insights into how importing and exporting countries can negotiate hydrogen trade agreements by structuring import order amounts based on country-specific characteristics. 
Finally, our sensitivity analyses identify the conditions required to achieve specific local production levels under the optimal dual sourcing policies. These conditions relate to production and storage costs, as well as variability in supply capacity, demand, and random yield. By linking the local production levels to the Netherlands' climate scenarios, which specify targets for local supply and import, our findings help policymakers assess the feasibility of meeting these targets and guide green hydrogen investment planning.

The remainder of the paper is structured as follows. Section~\ref{sec:literature} reviews the relevant literature. Section~\ref{sec:problem description} defines a  MDP model for the dual sourcing of hydrogen. Section~\ref{sec:solution approach} introduces heuristic dual sourcing policies. Section~\ref{sec:case study} details the case study focused on the Netherlands. Section~\ref{sec:results} presents the results of the case study, demonstrating the performance of the heuristic dual sourcing policies and providing decision-making insights based on sensitivity analyses. Finally, Section~\ref{sec:concluions} concludes with our key findings and final remarks.

\section{Literature}\label{sec:literature}
The study of dual sourcing in inventory management begins with \cite{barankin1961delivery} exploring one-period problems and \cite{daniel1962delivery} extending it to multiple periods. These initial studies provide a basis for understanding optimal inventory policies under specific conditions such as zero and one-unit lead times. For consecutive lead times where the difference between the lead times of suppliers is one, \cite{fukuda1964optimal} develops optimal inventory policies using a base-stock model that keeps inventory levels at a predetermined target. For general (not necessarily consecutive) lead times, \cite{whittemore1977optimal} show that the optimal dual sourcing policy becomes more complex and must consider the state vector including net inventory and outstanding orders. 

Despite extensive research efforts, finding an optimal dual sourcing policy under general lead times remains challenging. Several heuristic approaches have been developed to manage this complexity. For example, the dual index policy \citep{veeraraghavan2008now} monitors two inventory positions: one based on on-hand inventory and orders arriving from the faster but more expensive supplier, and the other from the slower but cheaper supplier. It aims to maintain optimal order-up-to levels for both positions by placing orders each period to bring inventory back to these targets. The capped dual index policy \citep{sun2019robust} builds on this by limiting the order quantity from the slower supplier, ensuring that no order exceeds a predetermined limit. The Vector Base-Stock (VBS) policy \citep{sheopuri2010new} uses a base-stock approach for orders from the faster supplier and a vector of base-stock levels for orders from the slower supplier, aimed at balancing inventory targets across multiple periods. The best-weighted bounds policy \citep{hua2015structural} extends the VBS approach by using a weighted average of upper and lower bounds for orders from the slower supplier. Although these heuristics do not guarantee an optimal solution, they provide practical policies to control inventory levels. The tailored base-surge policy \citep{allon2010global,janakiraman2015analysis} 
tracks a single inventory position, placing a fixed order with the slower supplier every period and only ordering from the faster supplier if the inventory falls below a set threshold, replenishing it to the desired level. \cite{xin2018asymptotic} demonstrate that tailored base-surge policy achieves asymptotic optimality as the lead time of the slower supplier increases, while the lead time of the faster supplier remains unchanged.




The consideration of supply uncertainties such as stochastic supply capacity and random yield presents further challenges. Stochastic supply capacity refers to random fluctuations in a supplier's available capacity, while random yield refers to cases where the delivered quantity is a random proportion of the ordered amount. Despite their significance, there has been relatively limited research on addressing these types of uncertainties in the dual sourcing context. For example,  \cite{jakvsivc2018dual} consider a dual sourcing inventory problem with random supply capacity at the faster supplier, where the supply capacity is not observed prior to placing orders to both supply sources in a particular period. They develop a myopic heuristic while assuming immediate delivery from the faster supplier and a one-period lead time from the slower supplier. \cite{chen2019heuristic} consider the dual sourcing problem with general lead time difference and random supply capacity for both sources.  They develop an LP-greedy heuristic based on techniques developed by \cite{chen2018preservation} and \cite{feng2018supply} to transform an inventory problem with random supply capacity into a convex optimization problem. 

Among random yield studies, \cite{chen2013sourcing} study a price-dependent demand model considering two unreliable suppliers with random yield and negligible lead times. They show that a reorder point exists for each supplier so that an order is placed for almost every inventory level below the reorder point. They also characterize the conditions under which optimal order quantities are decreasing in the inventory level, and a strict reorder point would apply. \cite{ju2015approximate} extend this work for general lead times, and develop a heuristic based on the above-mentioned dual-index order-up-to policy. \cite{gong2014dynamic} and \cite{zhu2015analysis} study a special case of random yield in the form of all-or-nothing supply disruptions. However, stochastic supply capacity and random yield have not been simultaneously considered in dual sourcing contexts. In this paper, we combine these two factors to analyze how they influence dual sourcing policies.

Several studies have examined the techno-economic aspects of local production and import for hydrogen supply rather than providing optimal dual sourcing policies. For example, \cite{eckl2022techno} compares on-site hydrogen production in southern Germany and importing hydrogen from Portugal. Their analysis assesses the economic viability and the logistical challenges associated with each option, emphasizing cost, infrastructure requirements, and energy efficiency. \cite{shin2023comparative} studies a comparative life cycle greenhouse gas analysis for clean hydrogen pathways in South Korea, comparing the environmental impacts of domestic production versus overseas import. Within Europe's net-neutral energy goals, \cite{frischmuth2022hydrogen} examines hydrogen sourcing strategies by analyzing the costs of domestic production versus import from outside Europe. However, these studies do not address how to optimize dual sourcing policies between local production and import. Our paper fills this gap by proposing a framework that accounts for key reliability challenges, including fluctuating supply capacities, variability in delivered quantities, and differences in lead times.





\section{Problem Description}\label{sec:problem description}

We consider a dual sourcing inventory control problem for green hydrogen. The system is situated on an infinite discrete-time horizon $T$. The system faces stochastic green hydrogen demand $d_t$ in each period $ t\in T$. Demand can be satisfied from two distinct sources: local production and import. The local production of green hydrogen is subject to stochastic supply capacity due to variable renewable energy output used in production. We assume that these stochastic capacities, denoted as $K^l_{t}$, are only revealed after an order is placed. This local supplier has a lead time $\tau^l$. In contrast, imported green hydrogen arrives from countries with strong renewable energy production capabilities. An import order for green hydrogen has a lead time of $\tau^i$ periods, and we assume $\tau^i>\tau^l$. Import is subject to random yield due to potential energy losses during transportation and conversion, resulting in the successful delivery of only a portion of the ordered quantity. The exact proportion, denoted by $p_t$, remains unknown until the arrival of the order. The local supplier has a unit ordering cost $c^l$, while the unit cost of import is $c^i$. The remaining hydrogen is carried over to the next period at a unit holding cost $c^h$. Unmet hydrogen demand is lost and incurs a unit penalty cost of $c^p$. The objective is to determine an optimal dual sourcing policy that minimizes the expected cost per period over the planning horizon.


\subsection{Markov Decision Process Description }


We now formulate our problem as a MDP. Decision epochs coincide with each period $t \in T$. The sequence of events within each period $t$ is as follows:
\begin{enumerate}   
    \item Orders placed $\tau^l \geq 1$ periods ago from the local supplier and $\tau^i \geq 1$ periods ago for imports arrive, and the inventory is updated.
    \item Based on the inventory level and the outstanding orders, the system places an order from the local supplier and the import source. In case of a zero lead time, the order arrives directly. Otherwise, it will arrive at the beginning of the periods $t +\tau^l$ or $t + \tau^i$, depending on the type of the order.
    \item Stochastic demand realizes and the supplier fulfills this demand with inventory in storage. Unmet demand incurs a unit lost sales cost, while excess inventory is carried over to the next period with a holding cost. Then, the system transitions to the next decision epoch.
\end{enumerate}

The system's state $\boldsymbol S_t$ at the beginning of period $t$ describes the inventory in storage $s_t^0$, the pipeline inventory from the local supplier $q_u^l$ for all $t - \tau^l  \leq u \leq t-1$, and the pipeline inventory from the import source $q_u^i$ for all $t - \tau^i  \leq u \leq t-1$. Thus, formally, we write $\boldsymbol S_t = (s_t^0, \boldsymbol q^l_t, \boldsymbol q^i_t)$. 

First, the pipeline inventory resulting from orders made at periods $t - \tau^l$ and $t - \tau^i$ (in case lead times are positive) will arrive. 
The import order $q_{t-\tau^i}^i$ arriving in period $t$ is subject to a proportional random yield. Only a random proportion $p_t$ of the ordered amount is received upon arrival. Specifically, the quantity obtained is $p_t \times q_{t-\tau^i}^i$. Here, $p_t$ is a realization of the random yield~$\tilde{p}$. 

After the pipeline inventory has arrived, the system can order hydrogen from the local supplier, denoted by $\hat{q}^l_t$, and by import, denoted by $q_t^i$. The local supplier has a stochastic capacity $\tilde{K}^l$ that limits the order quantity, where $K^l_t$ is a realization of this stochastic capacity at time $t$. Thus, the actual amount delivered equals $q^l_t = \min\{K^l_t, \hat{q}^l_t\}$, which will arrive in period $t + \tau^l$. Note that it will arrive immediately if $\tau^l = 0$.  
The import order $q_t^i$ will arrive in period $t + \tau^i$. The cost of ordering from the local supplier and import is given by the order cost function~$C^O(q^l_t, q^i_t) = c^lq^l_t + c^iq^i_t$.

Then, after orders are placed, we observe the demand $d_t$. Here $d_t$ is a realization of stochastic demand $\tilde{d}$ in period $t$. The transition towards state $\boldsymbol{S}_{t+1}$ is then obtained as
\begin{align}
    s'^0_t = s^0_t + q_{t - \tau^l}^l + p_tq_{t - \tau^i}^i - d_t, \\
    s^0_{t+1} = \max\{0, s'^0_t\},
\end{align}
The associated inventory costs equal $C^I(s'^0_t)=  c^h\max\{0, s'^0_t\} + c^p \max\{-s'^0_t, 0\}$, where $c^h$ is the unit holding cost and $c^p$ is the unit lost sales penalty cost. Note that $\boldsymbol S_{t+1} = (s^0_{t+1}, \boldsymbol q_{t+1}^l, \boldsymbol q_{t+1}^i)$, where $\boldsymbol q^l_{t+1}$ and $\boldsymbol q^i_{t+1}$ are obtained from $\boldsymbol q^l_t$ and $\boldsymbol q^i_t$ by removing the pipeline order that arrived at the start of period $t$ and inserting the order placed in period $t$. 

We can express the value function $V(\boldsymbol{S}_t)$, which represents the expected cost of being in state $\boldsymbol{S}_t$, as follows:\begin{align}
V(\boldsymbol{S}_t)= \min_{\hat{q}^l_t, q^i_t} \mathbb E \left [C^O(q^l_t, q^i_t) +  C^I(s'^{0}_t) + V(\boldsymbol S_{t+1})\right ].  
\end{align}
The goal is to find an optimal policy $\pi^*$ that prescribes an optimal ordering decision $\hat{q}^l_t, q^i_t$ that minimizes the expected cost per period. Clearly, this is a stationary policy, not dependent on $t$. The size of the state space depends on the lead times of local and import suppliers. As these lead times increase, the state space expands exponentially, creating computational challenges 
to find optimal solutions for MDPs. Thus, in the following section, in addition to the exact solution approaches, we propose heuristic dual sourcing policies. 

\section{Optimal and Heuristic Dual Sourcing Policies}\label{sec:solution approach}
In this section, we define various dual sourcing policies. First, we present a value iteration algorithm to derive an optimal dual sourcing policy for our MDP model. However, implementing this optimal policy in practice may face challenges due to potential variability in order quantities, which complicates production planning for suppliers that require stability. Therefore, we propose heuristic policies that approximate the optimal dual sourcing policy by stabilizing order quantities, which is particularly appealing for practitioners in complex supply chain environments. In addition, these heuristic policies can be a solution for the computational issues arising from the exponential growth of the state space with increased lead times. 



\subsection{Optimal Dual Sourcing Policy}

The value iteration algorithm is a fundamental approach for solving MDPs \citep{puterman2014markov,goyal2023first}. In our study, we use this algorithm to compute optimal dual sourcing policies, yielding a dynamic, state-dependent solution that determines optimal order quantities, as detailed in Algorithm~\ref{alg:value_iteration}. As our optimal policy is independent of time period $t$, we drop the period index in our state notation in the following section for notational convenience.

\begin{algorithm}
\caption{Value Iteration Algorithm}\label{alg:value_iteration}
\begin{algorithmic}[1]
\STATE \textbf{Output:} Optimal dual sourcing policy $\pi^*$
\STATE Initialize $V^0(\boldsymbol{S})$ for each state $\boldsymbol{S}$
\STATE $k \gets 0$
\WHILE{$\max_{\boldsymbol{S}} |V^k(\boldsymbol{S}) - V^{k-1}(\boldsymbol{S})| - \min_{\boldsymbol{S}} |V^k(\boldsymbol{S}) - V^{k-1}(\boldsymbol{S})| > \epsilon$}
    \FORALL{$ \boldsymbol{S} $}
        \FORALL{$\hat{q}^l,q^i$}
            \STATE $Q(\boldsymbol{S},\hat q^l,q^i) \gets \mathbb E[ C^O( q^l, q^i) + C^I(s'^{0}) + V^k(\boldsymbol S' )] $
        \ENDFOR
        \STATE $k \gets k + 1$
        \STATE $V^k(\boldsymbol{S}) \gets \min_{\hat{q}^l, q^i} Q(\boldsymbol{S},\hat q^l,q^i)$
        \STATE $\pi^*(\boldsymbol{S}) \gets \arg \min_{\hat{q}^l, q^i} Q(\boldsymbol{S},\hat q^l,q^i)$
    \ENDFOR
\ENDWHILE
\end{algorithmic}
\end{algorithm}

Algorithm~\ref{alg:value_iteration} derives an optimal dual sourcing policy $\pi^*$ by iteratively refining the value function for each state until convergence. Initially, the value function $ V^0(\boldsymbol{S})$ is set for all states, and an iteration counter $k$ is initialized at zero. For each state $\boldsymbol{S}$, the algorithm evaluates all possible combinations of order quantities $\hat{q}^l$ and $q^i$, calculating the total cost \( Q(\boldsymbol{S}, \hat q^i, q^l) \), which consists of the immediate cost and the expected future cost. The value function $V^k(\boldsymbol{S})$ and the optimal policy $\pi^*(\boldsymbol{S})$ are updated based on the actions that yield the lowest cost for each state. Note that we denote with $\boldsymbol S'$ the transitioned state. This process continues until the stopping criterion falls below a predefined threshold $\epsilon$.

The optimal dual sourcing policy, while theoretically cost-effective by obtaining the solution with the least cost, may encounter practical challenges in real-world applications. One potential concern is the variability in order amounts. For countries supplying hydrogen to multiple regions, such fluctuations can complicate the allocation of production. This could result in mismatches between supply and demand, leading to either surpluses or shortages. This issue is worsened by positive lead times, where fluctuating order sizes cause delays in deliveries, making it even harder to manage stock and maintain a stable supply chain. Additionally, local hydrogen production may require stability due to the shared use of renewable energy by other sectors, such as electricity generation or industrial processes. As a result, more stable hydrogen orders could help ensure efficient energy allocation across these different uses. In the following, we propose various heuristic dual sourcing strategies focusing on providing stable order amounts.

\subsection{Heuristic Dual Sourcing Policies}

To address the potential variability in order quantities arising from the optimal policy, we propose heuristic dual sourcing policies that facilitate implementation while maintaining cost-effectiveness.  One potential solution is inspired by Quantity Flexibility contracts \citep{tsay1999quantity,lian2009analysis,li2021supply}, which are commonly used in supply chains to improve coordination between suppliers and buyers. These contracts allow adjustments in order quantities within predefined limits for each planning period. This mechanism helps mitigate the risk of over- or under-supply by smoothing out fluctuations. In what follows, we introduce these policies and explain the calculations for policy-specific parameters.


\begin{enumerate} 
\item \texttt{FOQ}: The Fixed Order Quantity (FOQ) policy sets fixed quantities to be ordered both from local production $\bar{q}^l$ and import  $\bar{q}^i$. This policy serves as a baseline to evaluate the trade-off between stability and cost efficiency. 

To implement this policy, we run a simulation to determine the values of $\bar{q}^l, \bar{q}^i$.  We evaluate each combination of possible order quantities over 100,000 simulated periods and calculate the long-run average to identify the related optimal values.

\item \texttt{FOQ+}: The improved Fixed Order Quantity (FOQ+) policy determines and import quantities subject to adaptive adjustments within predefined ranges. We allow a discrete set of order quantities around the fixed values $\bar{q}^l,\bar{q}^i$ from \texttt{FOQ}. Specifically, the local order quantity can vary between $\bar{q}^l_{\text{min}}$ and $\bar{q}^l_{\text{max}}$, while the import quantity can vary between $\bar{q}^i_{\text{min}}$ and $\bar{q}^i_{\text{max}}$. This introduces some flexibility, allowing the policy to respond to uncertainties in a controlled manner.

We set the bounds $\bar{q}^l_{\text{min}}$, $\bar{q}^l_{\text{max}}$, $\bar{q}^i_{\text{min}}$, and $\bar{q}^i_{\text{max}}$, based on the optimal fixed order quantities $\bar{q}^i$ and $\bar{q}^l$ from \texttt{FOQ}. We update the action space of our MDP model by constraining $\hat{q}^l_t$ and $q^i_t$ within corresponding bounds. We then apply Algorithm \ref{alg:value_iteration}.

\item\texttt{TBS}: The Tailored Base-Surge (TBS) policy follows from the dual sourcing literature \citep{allon2010global,janakiraman2015analysis}. It orders a fixed import quantity, denoted by $\bar{q}^i$, while local orders are placed only when the inventory level drops below a specified threshold, $\delta^l$. The idea is to utilize local production as a backup supply source, reducing the reliance on fluctuating local capacity and making the policy more responsive to demand shortfalls. 

To determine the optimal fixed import quantity, $\bar{q}^i$, and inventory threshold, $\delta^l$, we apply the same simulation approach as in \texttt{FOQ}. We simulate various combinations of import quantities and threshold levels over 100,000 periods to calculate the long-run average cost, identifying the optimal values for each parameter.

\item\texttt{TBS+}: The improved Tailored Base-Surge (TBS+) policy builds on \texttt{TBS}. \texttt{TBS+} utilizes local supply when inventory falls below a threshold $\delta^l$ but introduces flexibility by allowing adjustments to the import quantity within a predefined range. Specifically, the import quantity can vary between a lower bound, $\bar{q}^i_{\text{min}}$, and an upper bound, $\bar{q}^i_{\text{max}}$. This policy allows for more responsive import quantities while maintaining the backup role of local production to manage inventory risks.

To obtain the TBS+ policy parameters, we adjust the MDP model by setting the action for local orders to $\max(0, \delta^l - s^0_t)$. For import orders, we limit the action space within $\bar{q}^i_{\text{min}}$ and $\bar{q}^i_{\text{max}}$, which are set based on $\bar{q}^i$ from \texttt{TBS}. We then apply Algorithm~\ref{alg:value_iteration}.

\end{enumerate}

\section{Case Study Description}\label{sec:case study}

In this section, we focus on a case study for the dual sourcing of green hydrogen within the Netherlands, where the transition to a green hydrogen economy is supported as part of the project ``Hydrogen Energy Applications in Valley Environments for Northern Netherlands" (HEAVENN) \citep{Heavenn}. We rely on data gathered from expert interviews within the HEAVENN project, government plans, 
technical reports and recent scientific literature. In the following subsections, we outline the current state of the Netherlands' hydrogen transition, local production dynamics, and the role of global hydrogen import in that context. 

\subsection{Transition to Green Hydrogen in the Netherlands}

In an effort to transition towards renewable energy sources, the Dutch government has committed to investments in hydrogen infrastructure \citep{NetherlandsHydrogenStrategy2020}. However, strategic analyses indicate that domestic hydrogen production will not be sufficient to meet future demand, even with planned expansions in renewable capacity \citep{PBL2022ClimateEnergyOutlook}. The variability in renewable energy supply, driven by inconsistent weather patterns, introduces uncertainty into hydrogen production through electrolysis. At the same time, the available renewable energy is needed across other sectors such as for electricity generation or for industrial processes, further complicating its availability for hydrogen production. Moreover, the cost of hydrogen production is highly dependent on electricity prices, with regions possessing higher renewable energy potential expected to produce hydrogen at lower costs than the Netherlands \citep{perey2023international}. Given these challenges, it has become clear that domestic production alone will not be sufficient to meet future hydrogen demand. As a result, the Netherlands will need to rely on hydrogen import from countries with more stable and abundant renewable energy resources. These regions, which benefit from consistently high levels of renewable energy, are expected to produce hydrogen at a lower and more predictable cost. However, importing hydrogen introduces additional challenges, including longer delivery lead times and potential energy losses during transport and conversion processes, which must be factored into strategic planning for the hydrogen supply chain \citep{IRENA2022Global}.

The Netherlands developed the climate scenarios IP2024 and II3050 to explore possible future energy pathways and support policymakers in making decisions about investment and infrastructure requirements \citep{ouden2020klimaatneutrale,netbeheer2023scenarios}.
Three different scenarios for 2030—Climate Ambition, National Drivers, and International Ambition—were created in collaboration with stakeholders and reflect current energy policies, market trends, and climate goals. Each scenario represents a different approach to achieving the Dutch target of at least a 55\% reduction in CO2 emissions by 2030. In the Climate Ambition scenario, domestic hydrogen production, driven by renewable energy like offshore wind, is significant but supplemented by import to meet rising demand. The National Drivers scenario prioritizes self-sufficiency, aiming to maximize local production through extensive government support for electrification and hydrogen generation, with minimal reliance on import. In contrast, the International Ambition scenario focuses on global cooperation, with the Netherlands relying mainly on hydrogen import from regions with stable renewable resources while still maintaining some local supply capacity.



In this case study, we aim to explore optimal dual sourcing policies, examining the trade-offs between local production and import. Specifically, we will analyze how the cost dynamics of local production, import, and storage, along with uncertainties such as stochastic local supply capacity, variable demand, random import yields, and lead time differences, affect optimal ordering decisions. These analyses can help policymakers make informed decisions on investments, ensuring that sourcing policies align with the Netherlands's three scenarios. In addition, our study aims to provide guidance to policymakers for international hydrogen trade agreements, examining how different scenarios affect the balance between local production and import.

The hydrogen market is evolving, and its dynamics are likely to shift as global trade expands and countries develop their hydrogen supply systems. In this study, we select 2030 as the base year, as many large-scale hydrogen production and infrastructure projects are targeting that timeframe, and predictive data for this period is widely available, providing a solid foundation for analysis. We consider an infinite horizon weekly ordering system. 
While we initially define the case study with specific parameters, later experiments will explore variability in these parameters, allowing for a broader analysis of potential hydrogen sourcing policies under different future scenarios.



\subsection{Local Production Dynamics}
The Dutch government plans to install 4 GW of electrolyzer capacity by 2030 in line with the Climate Ambition scenario \citep{pbl2023climate}. 
If this 4 GW electrolyzer capacity operates continuously at full capacity for one week, it could produce approximately 20,000 tonnes of hydrogen. To model potential variability in local production, we assume that stochastic capacity $K^l_t$ for each period $t$ is realized from a discrete set $\{0, 2000, 4000, 6000, 8000, 10000, 12000, 14000, 16000, 18000, 20000\}$. In our base case, we assume the electrolyzer operates at a 50\% capacity rate on average, and accordingly, we set the mean of the distribution to 10,000 tonnes. In our experiments, we consider different standard deviation values leading to various discrete probability distributions.

By 2030, hydrogen is expected to be used in the Netherlands across various sectors, including industry, transportation, the built environment, and power generation. 
Different reports project varying hydrogen demand for 2030 \citep{detz2020hydrogen}. 
Projections generally range between 75 and 100 petajoules.  We used the midpoint of this range, which corresponds to a weekly average of approximately 14,000 tonnes of hydrogen, as the mean in our demand distribution. To account for variability, we model demand $d_t$ for each period $t$ considering a discrete distribution, with possible values ranging from approximately 50\% below to 50\% above the mean of 14,000 tonnes: $\{6000, 8000, 10000, 12000, 14000, 16000, 18000, 20000, 22000\}$. Later, we varied the standard deviation to test different levels of variability in our experiments.


In our model, we assume that orders can be placed in increments of 2000 tonnes from both local production and import. This discretization can be practically reasonable in the context of hydrogen supply chains where bulk orders are typically placed due to logistical constraints. In addition, it is necessary to numerically solve the MDP model with Algorithm~\ref{alg:value_iteration} and the proposed heuristic dual sourcing policies in Section~\ref{sec:solution approach}.

In the Netherlands, hydrogen is planned to be stored both underground and above-ground. Underground storage will use salt caverns, like those being developed in Zuidwending \citep{hystock2023}. 
Above-ground options can be compressed gas tanks and liquid hydrogen tanks. 
The cost of storing hydrogen in 2030 is expected to vary depending on the storage type, technological advancements, and other factors. Based on reports by \citet{burke2024hydrogen} and \cite{gaster2024realist}, the cost ranges per kilogram per day in 2030 are determined as: €0.2 to €0.6 for salt caverns, €1.1 to €3.5 for compressed gas, and €2.0 to €5.0 for liquid hydrogen.  For the base case, we use the average of the projected cost ranges, multiply it by 7 to account for the weekly planning horizon, and use the resulting values for the parameter $c^h$. Later in the experiments, we test alternative cost scenarios to analyze potential fluctuations. We assume a lost sales penalty cost $c^p$ of €30 per kg of hydrogen, which is approximately four times the average unit ordering cost in exporting countries.

\subsection{Hydrogen Import and Global Trade Dynamics}

By 2030, the Netherlands plans to import hydrogen from several countries. These initiatives are supported by formal agreements, such as memorandums of understanding and joint statements, which outline collaborative efforts in hydrogen production and trade \citep{dutch_government}. In our case study, we focus on key countries such as Norway, Morocco, and the United Arab Emirates (UAE), not only for their strategic importance but also for their distinct characteristics in terms of lead times, production capacities, and transportation costs, which are crucial for understanding hydrogen import dynamics.

Hydrogen can be transported to the Netherlands in various ways, including pipelines, shipping, or even road transport, depending on the distance and infrastructure available. Pipeline transport via the European Hydrogen Backbone \citep{europeanhydrogenbackbone2020} is a potential option, as this network aims to connect hydrogen production and consumption sites across Europe. However, the full establishment of this infrastructure is not expected by 2030. As a result, in this case study, we focus solely on hydrogen transport by ship. Among the various shipping options—such as compressed gas or liquid hydrogen—we focus on ammonia, as its higher energy density makes it a potentially more efficient solution \citep{cha2018ammonia}.

The parameter settings for the exporting countries are presented in Table~\ref{Import Countries Data}, detailing distances, lead times, and transportation costs. We assume that ships depart from the main ports of the selected countries and arrive at the Port of Rotterdam, which is expected to serve as the primary hub for unloading and distributing hydrogen in the Netherlands. We calculate the distances and lead times using the SeaRates tool, which provides accurate transit times based on recent shipping data from various companies \footnote{SeaRates tool is available at \url{www.searates.com}}. Using these distances, we estimate the transportation costs for hydrogen to the Port of Rotterdam, incorporating conversion, transport, and re-conversion expenses based on \citet{IEA2019Future} and \citet{IRENA2022Global}. While our analysis focuses on dual sourcing strategies, it does not extend to the costs of distributing hydrogen within the Netherlands, whether it is imported through the Port of Rotterdam or produced locally.  
While these countries are actively investing in infrastructure to enhance their hydrogen production capabilities, the exact cost projections for 2030 are still evolving. For the base case, we set unit production cost per kilogram of hydrogen based on  \citet{cheng2024competitive} and \citet{Aurora2024}. 
We set the import cost parameter $c^i$ to the total cost in Table~\ref{Import Countries Data}. Later in our experiments, we consider a range of cost scenarios for local production $c^l$ and import $c^i$ to analyze their effect on the dual sourcing policies.


\begin{table}[h]
\centering \caption{Parameter setting for exporting countries}
\label{Import Countries Data}
\begin{tabular}{lllll}
\toprule
 & & \multicolumn{3}{c}{Costs (€/kg)} \\
 Country  &   Lead Time  & Production  & Transportation  & Total      \\
  \midrule
 Norway &  1 week
 &  6.21  & 2.41 & 8.62 \\
  Morocco  & 2 weeks 
 & 3.20  &  2.56 & 5.76 \\
  UAE    &  3 weeks
 & 3.60 &  2.67 &  6.27 \\
 \midrule 
\end{tabular}
\end{table}

Based on the report of \citet{IRENA2022Global}, for each import order, we assume a mean loss of 17.5\% and define possible arrived quantities constrained by a maximum loss of 35\%. To model the random yield of import, we assume a truncated normal distribution to the expected losses, with the mean set at 17.5\% and the distribution truncated between 0\% and 35\%, ensuring that the arrived quantities remain within realistic bounds. Using Monte Carlo simulation, we sample from the continuous truncated normal distribution. These continuous outcomes are then mapped to the nearest discrete values, which are in increments of 2000, ensuring the results align with the discrete nature of the system.



\section{Computational Results}\label{sec:results}

We present computational results in three parts. First, we evaluate the value of the dual sourcing model with stochastic local supply capacity and random yield by comparing it to single-sourcing models and dual sourcing models that do not account for stochastic supply, random yield, or both. Second, we assess the performance of the proposed heuristic dual sourcing policies across different parameter settings. Lastly, we conduct sensitivity analyses on production and storage cost parameters, as well as the variability level of uncertain factors, including local supply capacity, random import yield, and stochastic demand.

The Netherlands' climate scenarios—Climate Ambition, National Drivers, and International Ambition—target different levels of local production and import, aiming to explore future energy pathways to guide policymakers on investment and infrastructure planning. To provide insights aligned with these scenarios, we identify the necessary conditions to maintain specific levels of local production and import. 
Therefore, for the upcoming experiments, we report the percentage of local supply under various conditions. In addition, instead of using a fixed local production cost, we set it relative to the import cost through different ratios. We define the cost ratio $\rho_{l/i}$ as the ratio between the local production cost $c^l$ and the import cost $c^i$, given by $\rho_{l/i} = \frac{c^l}{c^i}$.

We introduce a variability level parameter, denoted $VarL$. The standard deviation for each distribution is defined as $ \sigma = VarL \times (\mu - L) $, where $\mu$ represents the mean and $L$ is the lower bound of the random variable, assuming the distributions are symmetric. For the base case, we assume stochastic supply capacity, demand, and random yield follow normal distributions with $VarL$ set to 0.5 while their means follow values reported in Section \ref{sec:case study}. 


\subsection{Value of dual sourcing with stochastic supply and random yield}

To demonstrate the value of our model, we compare it against several benchmark cases. First, we examine single sourcing cases where the system uses only the local source or only import. Additionally, to highlight the importance of accounting for both stochastic supply and random yield, we analyze three dual sourcing cases: (1) a case where both stochastic supply and random yield are disregarded, (2) a case where only random yield is ignored, and (3) a case where only stochastic supply is ignored. For each case, we solve the corresponding MDP models to obtain their optimal sourcing policy. Then, we simulate the system with stochastic supply and random yield under the obtained policies to evaluate the resulting cost.

We report the percentage cost deviation of each case relative to the optimal dual sourcing policies under both stochastic supply and random yield. Tables~\ref{tab:Norway_singlesourcing}, \ref{tab:Morocco_singlesourcing}, and \ref{tab:UAE_singlesourcing} present the results for various cost ratios and storage types for Norway, Morocco, and UAE, respectively. Note that differences in holding costs distinguish storage types, with salt caverns being the least expensive and liquid hydrogen the most. Abbreviations SC, CG, and LH refer to salt cavern, compressed gas, and liquid hydrogen storage, respectively. All cost values for the optimal policy (i.e., the last column in Tables~\ref{tab:Norway_singlesourcing}, \ref{tab:Morocco_singlesourcing}, and \ref{tab:UAE_singlesourcing}) are reported in thousands of euros per week.

\begin{table}[h]
\small
    \centering
     \caption{Percentage cost deviations relative to the optimal dual sourcing policy for Norway }\label{tab:Norway_singlesourcing}
    \resizebox{\columnwidth}{!}{%
    \begin{tabular}{ll|>{\raggedleft\arraybackslash}p{2.4cm}>{\raggedleft\arraybackslash}p{2.6cm}>{\raggedleft\arraybackslash}p{1.8cm}>{\raggedleft\arraybackslash}p{2cm}>{\raggedleft\arraybackslash}p{2cm}>{\raggedleft\arraybackslash}p{2cm}}
    \toprule
      &  &  \multicolumn{2}{c}{Single sourcing}  & \multicolumn{4}{c}{Dual sourcing: Stochastic supply/Random yield} \\
      \cmidrule(lr){3-4}  \cmidrule(lr){5-8} 
   $\rho_{l/i}$ & Storage  & Only local  (\%)& Only import  (\%) & No/No  (\%)  & Yes/No  (\%)  & No/Yes (\%) & Yes/Yes (€) \\ \midrule
\multirow{3}{*}{0.6} & SC & 76.54 & 53.68 & 17.89 & 2.79 & 13.30 & 97.71 \\
                     & CG & 44.64 & 36.28 & 13.27 & 2.98 & 7.86 & 121.51 \\
                     & LH & 37.19 & 33.14 & 10.72 & 4.43 & 11.38 & 129.53 \\
\midrule
\multirow{3}{*}{0.8} & SC & 65.08 & 30.71 & 12.45 & 2.41 & 8.55 & 114.93  \\
                     & CG & 41.31 & 21.26 & 7.53 & 4.90 & 5.44 & 136.56 \\
                     & LH & 35.56 & 19.91 & 8.75 & 3.64 & 3.69 & 143.82  \\
\midrule
\multirow{3}{*}{1.0} & SC & 57.49 & 14.34 & 7.38 & 9.21 & 5.13 & 131.41 \\
                     & CG & 39.71 & 10.08 & 6.04 & 3.70 & 4.25 & 150.45 \\
                     & LH & 34.81 & 9.58  & 3.73 & 2.23 & 3.16 & 157.41 \\
\midrule
\multirow{3}{*}{1.2} & SC & 53.39 & 2.76  & 3.03 & 3.03 & 2.40 & 146.12  \\
                     & CG & 41.63 & 3.16  & 0.70 & 0.62 & 4.33 & 160.44 \\
                     & LH & 37.80 & 3.61  & 1.20 & 1.20 & 5.14 & 166.40 \\
\midrule
\multirow{3}{*}{1.4} & SC & 61.23 & 0.38  & 1.01 & 1.01 & 0.00 & 149.51  \\
                     & CG & 49.38 & 1.13  & 0.91 & 0.91 & 0.18 & 163.69 \\
                     & LH & 44.92 & 1.42  & 0.32 & 0.32 & 0.06 & 170.08 \\
\midrule
\multicolumn{2}{c|}{Average} &  48.05 &	16.10 &	6.33 &	2.89 &	4.99 & 142.64  \\ \bottomrule
    \end{tabular}}
\end{table}

\begin{table}[h]
\small
    \centering
     \caption{Percentage cost deviations relative to the optimal dual sourcing policy for Morocco}\label{tab:Morocco_singlesourcing}
    \resizebox{\columnwidth}{!}{%
    \begin{tabular}{ll|>{\raggedleft\arraybackslash}p{2.4cm}>{\raggedleft\arraybackslash}p{2.6cm}>{\raggedleft\arraybackslash}p{1.8cm}>{\raggedleft\arraybackslash}p{2cm}>{\raggedleft\arraybackslash}p{2cm}>{\raggedleft\arraybackslash}p{2cm}}
    \toprule
      &  &  \multicolumn{2}{c}{Single sourcing}  & \multicolumn{4}{c}{Dual sourcing: Stochastic supply/Random yield} \\
      \cmidrule(lr){3-4}  \cmidrule(lr){5-8} 
   $\rho_{l/i}$ & Storage  & Only local  (\%)& Only import  (\%) & No/No  (\%)  & Yes/No  (\%)  & No/Yes (\%) & Yes/Yes (€) \\ \midrule
\multirow{3}{*}{0.6} & SC & 122.53 & 48.00 & 27.04 & 2.84 & 18.86 & 69.86 \\
                     & CG & 70.66 & 31.19 & 19.06 & 4.82 & 11.61 & 93.00 \\
                     & LH & 59.61 & 29.17 & 13.72 & 2.41 & 9.69 & 100.69 \\
\midrule
\multirow{3}{*}{0.8} & SC & 105.44 & 27.24 & 20.75 & 3.27 & 13.72 & 81.27 \\
                     & CG & 66.83 & 19.51 & 12.82 & 1.34 & 10.54 & 102.03 \\
                     & LH & 58.23 & 19.51 & 14.95 & 0.94 & 8.43 & 108.88 \\
\midrule
\multirow{3}{*}{1.0} & SC & 94.91 & 12.95 & 6.10 & 6.54 & 8.03 & 91.58 \\
                     & CG & 66.60 & 11.81 & 5.63 & 1.34 & 8.63 & 109.08 \\
                     & LH & 59.06 & 12.63 & 5.44 & 0.36 & 7.31 & 115.58 \\
\midrule
\multirow{3}{*}{1.2} & SC & 89.64 & 3.25 & 1.59 & 1.42 & 4.79 & 100.09 \\
                     & CG & 69.08 & 6.74 & 1.25 & 0.68 & 6.19 & 114.24 \\
                     & LH & 62.14 & 8.04 & 0.42 & 0.33 & 5.21 & 120.45 \\
\midrule
\multirow{3}{*}{1.4} & SC & 97.11 & 1.20 & 1.86 & 1.84 & 0.05 & 102.09 \\
                     & CG & 75.57 & 4.61 & 0.99 & 0.91 & 0.40 & 116.48 \\
                     & LH & 67.57 & 5.45 & 1.59 & 1.48 & 0.56 & 123.35 \\
\midrule
\multicolumn{2}{c|}{Average} & 77.67 &	16.09 &	8.87 &	2.04 &	7.60 & 103.24  \\ \bottomrule
    \end{tabular}}
\end{table}

\begin{table}[h]
\small
    \centering
    \caption{Percentage cost deviations relative to the optimal dual sourcing policy for UAE}\label{tab:UAE_singlesourcing}
        \resizebox{\columnwidth}{!}{%
        \begin{tabular}{ll|>{\raggedleft\arraybackslash}p{2.4cm}>{\raggedleft\arraybackslash}p{2.6cm}>{\raggedleft\arraybackslash}p{1.8cm}>{\raggedleft\arraybackslash}p{2cm}>{\raggedleft\arraybackslash}p{2cm}>{\raggedleft\arraybackslash}p{2cm}}
    \toprule
      &  &  \multicolumn{2}{c}{Single sourcing}  & \multicolumn{4}{c}{Dual sourcing: Stochastic supply/Random yield} \\
      \cmidrule(lr){3-4}  \cmidrule(lr){5-8} 
   $\rho_{l/i}$ & Storage  & Only local  (\%)& Only import  (\%) & No/No  (\%)  & Yes/No  (\%)  & No/Yes (\%) & Yes/Yes (€) \\ \midrule
\multirow{3}{*}{0.6} & SC & 108.36 & 48.33 & 28.93 & 3.28 & 20.02 & 76.06  \\
                     & CG & 62.63 & 33.23 & 16.74 & 4.22 & 12.75 & 99.14 \\
                     & LH & 52.70 & 31.23 & 12.83 & 1.90 & 15.18 & 107.35 \\
\midrule
\multirow{3}{*}{0.8} & SC & 93.65 & 27.73 & 22.65 & 4.50 & 14.95 & 88.34 \\
                     & CG & 59.34 & 21.16 & 14.93 & 2.37  &9.43 & 108.94  \\
                     & LH & 52.57 & 21.82 & 13.57 & 1.31 & 9.77 & 115.69  \\
\midrule
\multirow{3}{*}{1.0} & SC & 84.33 & 13.30 & 5.49 & 6.19 & 9.01 & 99.47  \\
                     & CG & 59.86 & 13.39 & 5.03 & 1.00 & 9.48 & 116.90  \\
                     & LH & 53.10 &	14.15 &	4.76 &	0.31 & 9.16 & 123.43	 \\
\midrule
\multirow{3}{*}{1.2} & SC & 80.74 & 4.02 & 2.04 & 1.81 & 6.13 & 108.49  \\
                     & CG & 61.76 & 7.49 & 0.75 & 0.73 & 7.12  & 122.73  \\
                     & LH & 56.00 & 9.07 & 0.21 & 0.13 & 6.52 & 128.89  \\
\midrule
\multirow{3}{*}{1.4} & SC & 87.46 & 1.37 & 1.49 & 1.49 &0.07 & 111.05  \\
                     & CG & 67.93 & 4.97 & 0.95 & 0.94 &0.04 &  125.72  \\
                     & LH & 61.55 & 6.34 & 1.20 & 1.17 &0.14 & 132.41  \\
\midrule
\multicolumn{2}{c|}{Average} &  69.47 &	17.17	& 8.77 &	2.09 & 8.65 & 110.97	 \\ \bottomrule
    \end{tabular}}
\end{table}

The ``Only local" and ``Only import" cases exhibit the largest cost deviations from the optimal dual sourcing policy across most cost ratios and storage types, underscoring the inefficiency of relying exclusively on a single source. The ``Only local" case incurs even higher deviations. At higher cost ratios ($\rho_{l/i} = 1.2$ and $\rho_{l/i} = 1.4$), where local production is more expensive than import, relying solely on local supply becomes inefficient, as expected. Conversely, at lower cost ratios, where local supply is cheaper, the "Only local" policy could be anticipated to perform better. However, since the local supplier has limited capacity, it cannot meet all demand, leading to lost sales costs for unmet demand, which drives up the cost deviation. 
On the other hand, the ``Only import" case achieves relatively lower deviations. There are cases, such as when $\rho_{l/i} = 1.2$ and $\rho_{l/i} = 1.4$, where deviations remain under 10\%, reflecting a relatively more cost-effective approach under these conditions.

When both stochastic supply and random yield are ignored (column No/No), the percentage deviations tend to decrease compared to single sourcing but are still considerably high, with values above $20\%$ such as for cost ratio of $\rho_{l/i} = 0.8$ and SC storage in Table~\ref{tab:Morocco_singlesourcing}. When random yield is ignored but stochastic supply is considered (column Yes/No), the deviations decrease further, particularly for lower cost ratios. For example, at $\rho_{l/i} = 0.6$ and $0.8$, the deviations are below 5\% in Tables~\ref{tab:Norway_singlesourcing} ,\ref{tab:Morocco_singlesourcing}, and \ref{tab:UAE_singlesourcing}.

When only random yield is considered (column No/Yes), the overall performance improves compared to the case where both stochastic supply and random yield are ignored. However, in some cases, such as at cost ratios of $\rho_{l/i} = 1$ and $\rho_{l/i} = 1.2$ (Tables~\ref{tab:Norway_singlesourcing}, \ref{tab:Morocco_singlesourcing}, and \ref{tab:UAE_singlesourcing}), ignoring these factors yields better outcomes.  This occurs because considering random yield reduces import reliability, while ignoring stochastic supply makes local supply appear more stable, shifting reliance toward local supply at these cost ratios. Therefore, the case that ignores reliability issues for both sources provides a more balanced comparison. In contrast, at a higher cost ratio of $\rho_{l/i} = 1.4$, the cost advantage of importing outweighs reliability concerns. As a result, the system primarily imports hydrogen, and the case with only random yield considered again outperforms the one that ignores both.

Overall, our case study demonstrates the advantages of dual sourcing over single sourcing in managing hydrogen supply. Considering both stochastic supply and random yield contributes to an average cost reduction of 8\% across all countries, with stochastic supply having a greater impact in this case study. However, the relative importance of these factors may vary for different case studies.


\subsection{Performance of the Heuristic Dual Sourcing Policies}\label{sec:policy comparison}

This section aims to identify heuristic dual sourcing policies that closely align with optimal policies by evaluating their performance under varying cost dynamics across countries. These insights can help decision-makers formulate hydrogen trade agreements with exporting countries by identifying the most effective policies for specific conditions. To construct the action set for local supply for \texttt{FOQ+}, we first obtain $\bar{q}^l$ from \texttt{TBS}. We then set $\bar{q}^l_{\text{min}}$ and $\bar{q}^l_{\text{max}}$ as the values located two steps to the left and right of $\bar{q}^l$ within the original action set. We follow the same procedure to construct the action sets for import for \texttt{FOQ+} and \texttt{TBS+}. For each dual sourcing policy, we present optimality gaps relative to the long-run average cost of optimal dual sourcing policies, along with the percentage of local supply to total supply. Tables~\ref{tab:Norway_policy}, \ref{tab:Morocco_policy}, and \ref{tab:UAE_policy} present the results for Norway, Morocco, and UAE, respectively. Note that all cost values for optimal policies are reported in thousands of euros per week.

\begin{table}[h]
\small
    \centering
    \caption{Performance of the policies for Norway}\label{tab:Norway_policy}
         \resizebox{\columnwidth}{!}{%
        \begin{tabular}{ll|rrrrrrrrrr}
    \toprule
    &  &  \multicolumn{2}{c}{\texttt{FOQ}}  & \multicolumn{2}{c}{\texttt{FOQ+}}  & \multicolumn{2}{c}{\texttt{TBS}}  & \multicolumn{2}{c}{\texttt{TBS+}} & \multicolumn{2}{c}{Optimal}    \\ 
    \cmidrule(lr){3-4} \cmidrule(lr){5-6} \cmidrule(lr){7-8} \cmidrule(lr){9-10} \cmidrule(lr){11-12}
   $\rho_{l/i}$ & Storage  & Gap (\%) &  Local (\%) & Gap (\%) &  Local (\%)   & Gap (\%) &  Local (\%)  & Gap (\%) &  Local (\%)   & Cost (€) & Local (\%)  \\ \midrule
    \multirow{3}{*}{0.6} & SC & 9.53 & 70.42 & 1.81 & 70.21 & 24.14 & 68.35 & 22.12 & 59.86 & 97.71 & 70.58 \\
                         & CG & 9.34 & 69.21 & 0.82 & 65.93 & 5.35 & 68.35 & 4.56 & 63.06 & 121.51 & 67.14 \\
                         & LH & 10.70 & 67.10 & 0.80 & 65.65 & 2.09 & 68.35 & 1.67 & 63.06 & 129.53 & 60.92 \\
    \midrule
    \multirow{3}{*}{0.8} & SC & 7.32 & 70.42 & 1.47 & 70.15 & 16.06 & 34.72 & 15.94 & 37.14 & 114.93 & 70.46 \\
                         & CG & 8.63 & 69.21 & 1.09 & 65.40 & 4.64 & 68.35 & 3.25 & 59.87 & 136.56 & 59.24 \\
                         & LH & 9.48 & 67.10 & 0.95 & 62.12 & 2.31 & 68.35 & 1.45 & 59.86 & 143.82 & 60.92 \\
    \midrule
    \multirow{3}{*}{1.0} & SC & 11.14 & 69.21 & 1.55 & 65.59 & 8.52 & 34.74 & 8.51 & 34.86 & 131.41 & 63.29 \\
                         & CG & 8.90 & 69.21 & 1.68 & 62.23 & 1.63 & 34.74 & 1.61 & 35.44 & 150.45 & 54.24 \\
                         & LH & 8.96 & 67.10 & 1.36 & 59.21 & 0.96 & 33.76 & 0.95 & 33.88 & 157.41 & 51.97 \\
    \midrule
    \multirow{3}{*}{1.2} & SC & 10.55 & 69.21 & 2.19 & 63.75 & 2.67 & 14.20 & 1.81 & 10.09 & 146.12 & 46.51 \\
                         & CG & 11.78 & 69.21 & 4.39 & 57.73 & 1.02 & 34.73 & 0.79 & 28.14 & 160.44 & 14.43 \\
                         & LH & 11.52 & 67.10 & 4.20 & 55.07 & 0.77 & 33.75 & 0.45 & 22.19 & 166.40 & 14.42 \\
    \midrule
    \multirow{3}{*}{1.4} & SC & 8.25 & 0.00 & 0.68 & 2.17 & 1.83 & 4.00 & 0.67 & 2.29 & 149.51 & 3.84 \\
                         & CG & 9.14 & 0.00 & 0.48 & 6.38 & 1.60 & 13.40 & 0.12 & 9.39 & 163.69 & 9.48 \\
                         & LH & 13.97 & 27.65 & 1.00 & 15.57 & 1.67 & 13.45 & 0.24 & 10.18 & 170.08 & 6.67 \\
    \midrule
       \multicolumn{2}{c|}{Average}   & 9.95 &	56.81 &	1.63 &	52.48 &	5.02 & 39.55 &	4.28 &	35.29 &	142.64	& 43.61 \\
    \bottomrule
    \end{tabular}}
\end{table}

\begin{table}[h]
\small
    \centering
     \caption{Performance of the policies for Morocco}\label{tab:Morocco_policy}
    \resizebox{\columnwidth}{!}{%
    \begin{tabular}{ll|rrrrrrrrrr}
    \toprule
    &  &  \multicolumn{2}{c}{\texttt{FOQ}}  & \multicolumn{2}{c}{\texttt{FOQ+}}  & \multicolumn{2}{c}{\texttt{TBS}}  & \multicolumn{2}{c}{\texttt{TBS+}} & \multicolumn{2}{c}{Optimal}    \\ 
    \cmidrule(lr){3-4} \cmidrule(lr){5-6} \cmidrule(lr){7-8} \cmidrule(lr){9-10} \cmidrule(lr){11-12} 
   $\rho_{l/i}$ & Storage  & Gap (\%) &  Local (\%) & Gap (\%) &  Local (\%)   & Gap (\%) &  Local (\%)  & Gap (\%) &  Local (\%)   & Cost (€) & Local (\%)  \\ \midrule
\multirow{3}{*}{0.6} & SC & 13.70 & 70.42 & 1.81 & 69.88 & 23.85 & 34.75 & 23.84 & 34.81 & 69.86 & 70.31 \\
                     & CG & 13.71 & 69.21 & 1.70 & 62.83 & 4.09 & 34.72 & 3.95 & 39.03 & 93.00 & 57.54 \\
                     & LH & 17.15 & 67.10 & 1.58 & 57.38 & 2.36 & 46.19 & 0.54 & 49.20 & 100.69 & 52.26 \\
\midrule
\multirow{3}{*}{0.8} & SC & 11.14 & 70.42 & 1.51 & 68.10 & 13.98 & 34.73 & 13.98 & 34.73 & 81.27 & 67.29 \\
                     & CG & 13.80 & 69.21 & 2.31 & 59.57 & 0.91 & 34.74 & 0.91 & 34.74 & 102.03 & 49.51 \\
                     & LH & 16.96 & 67.10 & 2.38 & 54.21 & 0.55 & 34.75 & 0.54 & 34.81 & 108.88 & 45.57 \\
\midrule
\multirow{3}{*}{1}   & SC & 10.53 & 70.42 & 1.93 & 64.92 & 7.87 & 34.74 & 7.68 & 25.25 & 91.58 & 60.82 \\
                     & CG & 15.93 & 69.21 & 4.33 & 55.14 & 0.01 & 34.74 & 0.01 & 34.74 & 109.08 & 34.90 \\
                     & LH & 18.31 & 67.10 & 4.00 & 50.86 & 0.10 & 33.74 & 0.10 & 33.74 & 115.58 & 34.90 \\
\midrule
\multirow{3}{*}{1.2} & SC & 13.16 & 12.27 & 0.97 & 11.52 & 1.99 & 14.23 & 1.10 & 8.20 & 100.09 & 26.23 \\
                     & CG & 14.64 & 0.00 & 1.45 & 8.19 & 0.58 & 24.63 & 0.21 & 21.39 & 114.24 & 19.73 \\
                     & LH & 20.82 & 27.65 & 0.10 & 18.24 & 0.26 & 23.68 & 0.03 & 20.26 & 120.45 & 18.73 \\
\midrule
\multirow{3}{*}{1.4} & SC & 14.05 & 0.00 & 0.60 & 3.17 & 1.22 & 4.32 & 0.47 & 3.64 & 102.09 & 4.16 \\
                     & CG & 12.44 & 0.00 & 0.79 & 7.41 & 0.97 & 13.41 & 0.11 & 9.39 & 116.48 & 9.44 \\
                     & LH & 13.68 & 0.00 & 0.79 & 9.75 & 0.63 & 13.42 & 0.03 & 12.65 & 123.35 & 12.93 \\
\midrule
       \multicolumn{2}{c|}{Average}   & 14.67 &44.01 &	1.75 & 40.08 &	3.96 &	27.79 &	3.57 &	26.44 &	103.24 & 37.62 \\ \bottomrule
    \end{tabular}}
\end{table}

\begin{table}[h]
\small
    \centering
    \caption{Performance of the policies for UAE}\label{tab:UAE_policy}
    \resizebox{\columnwidth}{!}{%
    \begin{tabular}{ll|rrrrrrrrrr}
    \toprule
    &  &  \multicolumn{2}{c}{\texttt{FOQ}}  & \multicolumn{2}{c}{\texttt{FOQ+}}  & \multicolumn{2}{c}{\texttt{TBS}}  & \multicolumn{2}{c}{\texttt{TBS+}} & \multicolumn{2}{c}{Optimal}    \\ 
    \cmidrule(lr){3-4} \cmidrule(lr){5-6} \cmidrule(lr){7-8} \cmidrule(lr){9-10} \cmidrule(lr){11-12}
   $\rho_{l/i}$ & Storage  & Gap (\%) &  Local (\%) & Gap (\%) &  Local (\%)   & Gap (\%) &  Local (\%)  & Gap (\%) &  Local (\%)   & Cost (€) & Local (\%)  \\ \midrule
    \multirow{3}{*}{0.6} & SC & 10.99 & 70.42 & 1.40 & 70.09 & 22.60 & 34.74 & 22.61 & 34.80 & 76.06 & 70.63 \\
                         & CG & 11.24 & 69.21 & 1.46 & 63.45 & 4.43 & 46.20 & 2.41 & 50.43 & 99.14 & 58.12 \\
                         & LH & 13.71  & 67.10 & 1.58  & 57.10  & 2.29 & 46.20 & 0.97 & 50.59  & 107.35  & 52.69 \\
    \midrule
    \multirow{3}{*}{0.8} & SC & 8.97 & 70.42 & 1.23 & 66.66 & 13.10 & 34.74 & 13.10 & 34.74 & 88.34 & 64.29 \\
                         & CG & 11.58 & 69.21 & 2.27 & 60.26 & 1.19 & 34.74 & 1.19 & 34.74 & 108.94 & 50.28 \\
                         & LH & 14.72 & 67.10 & 2.20 & 56.14 & 0.91 & 34.74 & 0.92 & 34.80 & 115.69 & 50.86 \\
    \midrule
    \multirow{3}{*}{1.0} & SC & 8.68 & 70.42 & 1.68 & 64.36 & 7.16 & 34.74 & 7.09 & 30.71 & 99.47 & 61.27 \\
                         & CG & 13.61 & 69.21 & 4.25 & 57.19 & 0.00 & 34.74 & 0.00 & 34.74 & 116.90 & 34.75 \\
                         & LH & 15.81 & 67.10 & 3.73 & 52.86 & 0.00 & 33.75 & 0.00 & 33.75 & 123.43 & 33.75 \\
    \midrule
    \multirow{3}{*}{1.2} & SC & 13.00 & 12.27 & 1.14 & 10.56 & 2.01 & 14.26 & 1.24 & 8.01 & 108.49 & 25.69 \\
                         & CG & 11.71 & 0.00 & 1.76 & 9.51 & 0.70 & 13.45 & 0.44 & 17.19 & 122.73 & 19.23 \\
                         & LH & 18.68 & 27.65 & 0.13 & 20.19 & 0.23 & 23.69 & 0.00 & 20.83 & 128.89 & 20.83 \\
    \midrule
    \multirow{3}{*}{1.4} & SC & 14.14 & 0.00 & 0.58 & 3.67 & 0.40 & 6.09 & 0.34 & 5.91 & 111.05 & 4.96 \\
                         & CG & 9.08 & 0.00 & 0.74 & 7.85 & 0.49 & 13.46 & 0.05 & 10.66 & 125.72 & 10.89 \\
                         & LH & 8.50 & 0.00 & 0.82 & 10.54 & 0.34 & 13.49 & 0.01 & 14.85 & 132.41 & 14.88 \\
    \midrule
       \multicolumn{2}{c|}{Average}   & 12.29 & 44.01 &	1.66 & 40.70 & 3.72 & 27.94	& 3.36 & 27.78 &	110.97 &	38.21  \\
    \bottomrule
    \end{tabular}}
\end{table}

Among the heuristic policies, \texttt{FOQ+} demonstrates the best overall performance, followed by \texttt{TBS+}, \texttt{TBS}, and finally \texttt{FOQ}. \texttt{FOQ}'s relatively weak performance is expected, as it always orders fixed amounts and cannot adapt to changing conditions. However, \texttt{FOQ} outperforms \texttt{TBS} in a few specific scenarios where both the cost ratio and storage costs are low, such as when $\rho_{l/i} = 0.6$ and SC storage is used in all countries. In such cases, the optimal policy increasingly favors the local supplier. While \texttt{TBS} relies on local supply as a backup source only when inventory falls below a threshold, \texttt{FOQ} takes a more balanced approach, resulting in better performance.

Incorporating a degree of flexibility into \texttt{FOQ} (resulting in \texttt{FOQ+}) substantially enhances its performance, yielding an average cost benefit of 11\% compared to \texttt{FOQ} across all countries. Specifically, \texttt{FOQ+} has an average gap of 1.63\%, 1.75\%, and 1.66\% in Tables~\ref{tab:Norway_policy}, \ref{tab:Morocco_policy}, and \ref{tab:UAE_policy}, respectively. This indicates its robustness across a range of scenarios, making \texttt{FOQ+} a strong candidate for trade agreements. Nevertheless, there are many scenarios where \texttt{FOQ+} is outperformed by \texttt{TBS+} and even by \texttt{TBS}. This mostly happens when $\rho_{l/i} \geq 1$ and storage costs are relatively high, as in the CG and LH cases. 

\texttt{TBS} and \texttt{TBS+} perform particularly well as the cost ratio and storage costs increase. For example, in Tables~\ref{tab:Morocco_policy} and \ref{tab:UAE_policy}, both policies achieve a gap below 1\% in more than 50\% of scenarios. However, their average performance is relatively worse for Norway, as shown in Table~\ref{tab:Norway_policy}. These results align with the findings in the literature \citep{xin2018asymptotic}, which confirm that the performance of TBS policies improves with increasing lead times. When comparing \texttt{TBS} and \texttt{TBS+}, we observe notable improvements. For example, in Table~\ref{tab:Norway_policy}, for $\rho_{l/i} = 1.4$ with CG storage, the gap decreases from 1.60\% to 0.12\%. The performance of \texttt{TBS} and \texttt{TBS+} declines as the cost ratio and storage cost decrease. This occurs because these policies treat the local source as a backup, resulting in a lower local supply percentage compared to the optimal policy, see $\rho_{l/i} = 0.6$ and $\rho_{l/i} = 0.8$ scenarios in all countries.

Overall, \texttt{FOQ+} demonstrates the best average performance, while \texttt{TBS+} achieves superior performance where the cost ratio and storage cost are high. These findings underscore the importance of tailoring policy selection to the cost dynamics and country characteristics in order to make effective hydrogen trade agreements between importing and exporting countries.


\subsection{Sensitivity Analyses}\label{sec:sensitivity_analyses}

For the sensitivity analyses, we focus on Morocco as the exporting country due to its low anticipated production costs and reasonable lead time, making it an economically favorable option for hydrogen import to the Netherlands. Under optimal dual sourcing policies, we analyze the effect of several parameters on the local supply rates. Specifically, Section~\ref{sec:production cost} focuses on the effect of production costs, while Section~\ref{sec:storage cost} focuses on the effect of storage costs. Finally, Section~\ref{sec:uncertainty} investigates the effect of variability in the stochastic parameters, including local supply capacity, random import yield, and stochastic demand.

\subsubsection{Production cost impact on local supply rate dynamics.} \label{sec:production cost}
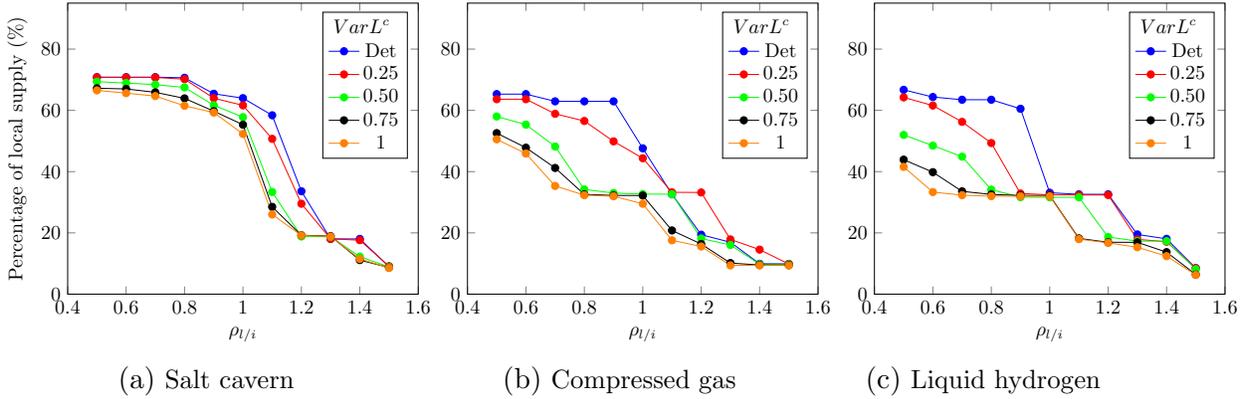
\begin{figure}[h]
\begin{subfigure}{0.33\textwidth}
    \centering
        \begin{tikzpicture}[scale=0.68, transform shape]
            \begin{axis}[
        xlabel={$\rho_{l/i}$},
        ylabel={Percentage of local supply (\%) },
        ymin=0, ymax=95,
        xticklabel style={/pgf/number format/1000 sep=},
        legend pos=north east,
        ymajorgrids=false,
        grid style=dashed,
    ]
    \addlegendimage{empty legend}
    \addplot[blue,mark=*] coordinates {
    (0.5, 70.84177533)
    (0.6, 70.84177533)
    (0.7, 70.83930241)
    (0.8, 70.62743104)
    (0.9, 65.41702797)
    (1.0, 63.98433917)
    (1.1, 58.39690682)
    (1.2, 33.59490764)
    (1.3, 18.02748307)
    (1.4, 18.02748307)
    (1.5, 8.962489303)
};

\addplot[red,mark=*] coordinates {
    (0.5, 70.7460875)
    (0.6, 70.7460875)
    (0.7, 70.7460875)
    (0.8, 70.1787755)
    (0.9, 63.98176435)
    (1.0, 61.6389071)
    (1.1, 50.69927019)
    (1.2, 29.51965142)
    (1.3, 18.05369416)
    (1.4, 17.61629049)
    (1.5, 9.133561213)
};

\addplot[green,mark=*] coordinates {
    (0.5, 69.36009052)
    (0.6, 68.92228338)
    (0.7, 68.40122923)
    (0.8, 67.48998592)
    (0.9, 61.67735752)
    (1.0, 57.79032248)
    (1.1, 33.3571765)
    (1.2, 18.82190004)
    (1.3, 18.71122199)
    (1.4, 12.25221592)
    (1.5, 8.885578145)
};

\addplot[black,mark=*] coordinates {
    (0.5, 67.21620748)
    (0.6, 66.99525769)
    (0.7, 65.89078449)
    (0.8, 63.91862033)
    (0.9, 59.73773707)
    (1.0, 55.31424098)
    (1.1, 28.4916619)
    (1.2, 19.23166644)
    (1.3, 18.95013478)
    (1.4, 11.12902543)
    (1.5, 8.673039893)
};

\addplot[orange,mark=*] coordinates {
    (0.5, 66.52742907)
    (0.6, 65.62467257)
    (0.7, 64.66971578)
    (0.8, 61.498885)
    (0.9, 59.28942483)
    (1.0, 52.29262077)
    (1.1, 26.01335359)
    (1.2, 19.17099789)
    (1.3, 18.80994349)
    (1.4, 11.49537421)
    (1.5, 8.575432636)
};
     \addlegendentry{\hspace{-.75cm}\textbf{$VarL^c$}}
     \addlegendentry{Det} 
      \addlegendentry{0.25} 
       \addlegendentry{0.50} 
        \addlegendentry{0.75} 
        \addlegendentry{1} 
          \end{axis}
        \end{tikzpicture}
    \caption{\small Salt cavern}\label{fig:DifferentLevelsSalt} 
\end{subfigure}
\begin{subfigure}{0.32\textwidth}
    \centering
        \begin{tikzpicture}[scale=0.68, transform shape]
            \begin{axis}[
        xlabel={$\rho_{l/i}$},
        ymin=0, ymax=95,
        xticklabel style={/pgf/number format/1000 sep=},
        legend pos=north east,
        ymajorgrids=false,
        grid style=dashed,
    ]
    \addlegendimage{empty legend}

   \addplot[blue,mark=*] coordinates {
    (0.5, 65.30589635)
    (0.6, 65.30589635)
    (0.7, 62.94733224)
    (0.8, 62.94733224)
    (0.9, 62.94733224)
    (1.0, 47.61115786)
    (1.1, 32.59356138)
    (1.2, 19.40371535)
    (1.3, 16.95050004)
    (1.4, 9.918460842)
    (1.5, 9.886932419)
};

\addplot[red,mark=*] coordinates {
    (0.5, 63.64564566)
    (0.6, 63.64564566)
    (0.7, 58.88367967)
    (0.8, 56.53558886)
    (0.9, 49.87674864)
    (1.0, 44.37481312)
    (1.1, 33.24538063)
    (1.2, 33.19733466)
    (1.3, 17.8722724)
    (1.4, 14.53546965)
    (1.5, 9.822515054)
};

\addplot[green,mark=*] coordinates {
    (0.5, 58.00729133)
    (0.6, 55.33330589)
    (0.7, 48.20013028)
    (0.8, 34.18256661)
    (0.9, 33.00951661)
    (1.0, 32.66273483)
    (1.1, 32.62066921)
    (1.2, 18.17528256)
    (1.3, 16.04801568)
    (1.4, 9.770592689)
    (1.5, 9.736606481)
};

\addplot[black,mark=*] coordinates {
    (0.5, 52.56923553)
    (0.6, 47.81029728)
    (0.7, 41.20004437)
    (0.8, 32.56544803)
    (0.9, 32.22622159)
    (1.0, 32.18582233)
    (1.1, 20.79206759)
    (1.2, 16.37092806)
    (1.3, 10.14836515)
    (1.4, 9.520709578)
    (1.5, 9.484988466)
};

\addplot[orange,mark=*] coordinates {
    (0.5, 50.55609574)
    (0.6, 45.90275787)
    (0.7, 35.32773007)
    (0.8, 32.31931623)
    (0.9, 31.95601743)
    (1.0, 29.51342756)
    (1.1, 17.57307042)
    (1.2, 15.55740288)
    (1.3, 9.393295218)
    (1.4, 9.389656606)
    (1.5, 9.358326584)
};
     \addlegendentry{\hspace{-.75cm}\textbf{$VarL^c$}}
     \addlegendentry{Det} 
      \addlegendentry{0.25} 
       \addlegendentry{0.50} 
        \addlegendentry{0.75} 
        \addlegendentry{1} 
          \end{axis}
        \end{tikzpicture}
    \caption{\small Compressed gas}\label{fig:DifferentLevelsGas}
\end{subfigure}
\begin{subfigure}{0.25\textwidth}
    \centering
        \begin{tikzpicture}[scale=0.68, transform shape]
            \begin{axis}[
        xlabel={$\rho_{l/i}$},
        ymin=0, ymax=95,
        xticklabel style={/pgf/number format/1000 sep=},
        legend pos=north east,
        ymajorgrids=false,
        grid style=dashed,
    ]
    \addlegendimage{empty legend}
    
\addplot[blue,mark=*] coordinates {
    (0.5, 66.70882198)
    (0.6, 64.35519738)
    (0.7, 63.47064092)
    (0.8, 63.47064092)
    (0.9, 60.52031102)
    (1.0, 33.10047781)
    (1.1, 32.59401081)
    (1.2, 32.58752678)
    (1.3, 19.48884442)
    (1.4, 18.04864241)
    (1.5, 8.503961365)
};

\addplot[red,mark=*] coordinates {
    (0.5, 64.24340482)
    (0.6, 61.58000525)
    (0.7, 56.24785063)
    (0.8, 49.33408072)
    (0.9, 32.8432855)
    (1.0, 32.42334055)
    (1.1, 32.37554313)
    (1.2, 32.3696946)
    (1.3, 17.87864651)
    (1.4, 17.11922742)
    (1.5, 8.461909912)
};

\addplot[green,mark=*] coordinates {
    (0.5, 51.98515675)
    (0.6, 48.47996275)
    (0.7, 44.88969668)
    (0.8, 34.12107439)
    (0.9, 31.65656593)
    (1.0, 31.61299344)
    (1.1, 31.60729729)
    (1.2, 18.64732501)
    (1.3, 17.31781576)
    (1.4, 17.30003212)
    (1.5, 8.129764468)
};

\addplot[black,mark=*] coordinates {
    (0.5, 43.90580152)
    (0.6, 39.82159322)
    (0.7, 33.55952527)
    (0.8, 32.5505528)
    (0.9, 32.22520002)
    (1.0, 32.18582233)
    (1.1, 18.20156015)
    (1.2, 16.93614914)
    (1.3, 16.91647693)
    (1.4, 13.71012773)
    (1.5, 6.403627948)
};

\addplot[orange,mark=*] coordinates {
    (0.5, 41.54517124)
    (0.6, 33.31788658)
    (0.7, 32.31016841)
    (0.8, 31.99523888)
    (0.9, 31.95601743)
    (1.0, 31.95405775)
    (1.1, 17.95828475)
    (1.2, 16.74352524)
    (1.3, 15.30014898)
    (1.4, 12.40113366)
    (1.5, 6.276405666)
};

     \addlegendentry{\hspace{-.75cm}\textbf{$VarL^c$}}
     \addlegendentry{Det} 
      \addlegendentry{0.25} 
       \addlegendentry{0.50} 
        \addlegendentry{0.75} 
        \addlegendentry{1} 
          \end{axis}
        \end{tikzpicture}
    \caption{\small Liquid hydrogen }\label{fig:DifferentLevelsLiquid}
\end{subfigure}
\caption{Impact of the hydrogen production price dynamics on the local supply rate for various levels of $VarL^c$}\label{fig:DifferentLevels}
\end{figure}

We conduct a series of experiments across different production cost ratios to assess the impact on the percentage of hydrogen sourced from local supplier (i.e., the local supply rate). Figure~\ref{fig:DifferentLevels} shows the local supply rate across different production cost ratios ($\rho_{l/i}$) and variability levels of stochastic capacity ($VarL^c$) for each storage type.

As expected, the percentage of local supply decreases as both the cost ratio ($\rho_{l/i}$) and the variability level of stochastic capacity ($VarL^c$) increase. When $\rho_{l/i}$ and $VarL^c$ are rather low, one might expect the percentage of local supply to approach 100\%. Yet, this is not observed due to limitations in electrolyzer capacity and capacity rate, as we assume a climate ambition scenario with a 4GW electrolyzer operating at a 50\% capacity rate on average. This capacity constraint limits the system's ability to fully rely on local production, even when costs are favorable. This indicates that if local hydrogen production costs can be reduced to certain levels, the system would have the tendency to increase its reliance on local production. In such a scenario, decision-makers would be encouraged to increase investments in electrolyzer capacities.

Interestingly, the decline in the percentage of local supply is not uniform; while some ranges show gradual decreases, others exhibit sharp drops, indicating critical points in the cost structure. These sharp changes typically occur around $\rho_{l/i}=1$, which reflects the point where local production and import costs become more comparable.

In all figures, at higher $\rho_{l/i}$, the system’s reliance on local production becomes low. 
However, local supply does not fully drop to zero, as the system still relies on the faster response time of local supply to handle stochastic demand, particularly in urgent situations.

These variations in the percentage of local supply, as shown in our graphs, are essential for evaluating the feasibility of the Netherlands' climate scenarios. For instance, to achieve the National Drivers scenario, which emphasizes a higher reliance on local hydrogen production, low cost ratios and minimal variability in stochastic capacity are required. Additionally, if hydrogen can be stored at low costs, such as in the case of salt caverns, the National Drivers scenario is more attainable. On the other hand, the International Ambition scenario, which relies more heavily on importing, becomes more feasible as $\rho_{l/i}$ and $VarL^c$ increase, especially when local production costs exceed import costs by 20\% or more, as seen around $\rho_{l/i} = 1.2$. Importantly, our analysis highlights that variations in cost ratios and capacity variability result in significant differences in the percentage of local supply. These insights are essential for guiding decision-makers on the conditions and investments required to effectively align with the ambitions of each climate scenario.

\subsubsection{Impact of storage cost on local supply rate.} \label{sec:storage cost}

\begin{figure}[h]
    \centering
        \begin{tikzpicture}[scale=0.68, transform shape]
            \begin{axis}[
        xlabel={Storage cost (€/kg)},
        ylabel={Percentage of local supply (\%) },
        ymin=0, ymax=105,
        xticklabel style={/pgf/number format/1000 sep=},
        legend pos=north east,
        ymajorgrids=false,
        grid style=dashed,
    ]
    \addlegendimage{empty legend}
    
\addplot[blue,mark=*] coordinates {
    (0.5, 69.40952117)
    (1, 63.28081534)
    (1.5, 57.29121613)
    (2, 55.33330589)
    (2.5, 52.73724739)
    (3, 49.399402)
    (3.5, 48.47996275)
    (4, 48.44821074)
    (4.5, 46.77960157)
    (5, 46.77960157)
    (5.5, 46.77960157)
    (6, 46.77960157)
    (6.5, 46.77960157)
    (7, 45.35747028)
    (7.5, 46.52949817)
};
\addplot[red,mark=*] coordinates {
    (0.5, 64.24889931)
    (1, 56.79240467)
    (1.5, 48.80993067)
    (2, 42.07485522)
    (2.5, 34.12107439)
    (3, 34.12107439)
    (3.5, 34.12107439)
    (4, 32.01360513)
    (4.5, 32.01360513)
    (5, 32.01360513)
    (5.5, 32.01360513)
    (6, 32.01360513)
    (6.5, 31.99283933)
    (7, 31.9968885)
    (7.5, 31.9968885)
    (8, 29.68285919)
};

\addplot[green,mark=*] coordinates {
    (0.5, 53.9871189)
    (1, 33.50193302)
    (1.5, 32.66344076)
    (2, 32.66344076)
    (2.5, 32.66273483)
    (3, 32.66273483)
    (3.5, 31.61299344)
    (4, 31.6128895)
    (4.5, 31.61334547)
    (5, 31.61334547)
    (5.5, 31.61334547)
    (6, 31.61334547)
    (6.5, 31.65292975)
    (7, 31.65292975)
    (7.5, 29.63607617)
    (8, 29.63607617)
};

\addplot[black,mark=*] coordinates {
    (0.5, 18.82190004)
    (1, 17.80446068)
    (1.5, 18.20442324)
    (2, 18.07798695)
    (2.5, 18.17528256)
    (3, 18.55583383)
    (3.5, 18.64732501)
    (4, 31.60603334)
    (4.5, 31.60729729)
    (5, 31.60729729)
    (5.5, 31.60729729)
    (6, 31.60729729)
    (6.5, 31.6110656)
    (7, 29.63068955)
    (7.5, 29.63068955)
    (8, 29.63068955)
};

\addplot[orange,mark=*] coordinates {
    (0.5, 11.67237701)
    (1, 9.749661039)
    (1.5, 9.749661039)
    (2, 9.770592689)
    (2.5, 10.02714792)
    (3, 15.00559131)
    (3.5, 17.30003212)
    (4, 17.31169669)
    (4.5, 15.56120419)
    (5, 15.60427413)
    (5.5, 16.6860673)
    (6, 15.44688355)
    (6.5, 15.45311494)
    (7, 15.45311494)
    (7.5, 16.67612178)
    (8, 16.75183753)
};

     \addlegendentry{\hspace{-.75cm}$\boldsymbol{\rho_{l/i}}$}
     \addlegendentry{0.6} 
      \addlegendentry{0.8} 
       \addlegendentry{1} 
        \addlegendentry{1.2} 
        \addlegendentry{1.4} 
          \end{axis}
        \end{tikzpicture}
    \caption{Impact of the storage cost on the local supply rate }\label{fig:DiffStorageCost}
\end{figure}
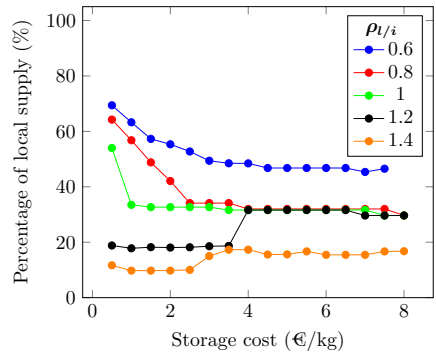

Figure~\ref{fig:DifferentLevels} shows that the type of storage also plays an important role. When storage costs are low, such as in the case of salt caverns, and local production remains relatively inexpensive, the system relies more heavily on local production compared to other storage technologies. To better understand the impact of storage costs, we present the results for a fixed $VarL^c$ of 0.5 under varying storage cost conditions in Figure~\ref{fig:DiffStorageCost}.

Figure~\ref{fig:DiffStorageCost} shows that when storage costs become high, particularly as they approach or exceed production costs, there are no significant changes in the percentage of local supply rate. This indicates that once storage costs exceed a certain threshold, further increases no longer impact the local supply rate. Specifically, for storage costs exceeding 4€/kg, we observe that the local supply rate remains quite similar when $\rho_{l/i}$ is 0.6, 0.8, and 1. At this point, storage costs dominate production cost dynamics, becoming the primary factor in the system's reliance on local production. 

However, when storage costs are relatively low, we observe significant changes in the optimal dual sourcing strategy. 
When local production costs are lower, increasing storage costs reduces reliance on local production. Conversely, when import costs are lower, increasing storage costs leads to greater reliance on local production. These trends show that as storage costs increase, the system reduces over-reliance on any single supply source to minimize the amount kept in storage. This is particularly crucial when both sources present disadvantages, such as stochastic supply, random yield, and positive lead time, since over-reliance on either source would necessitate holding larger amounts in storage.

\subsubsection{Analysis of the uncertain factors.}\label{sec:uncertainty}

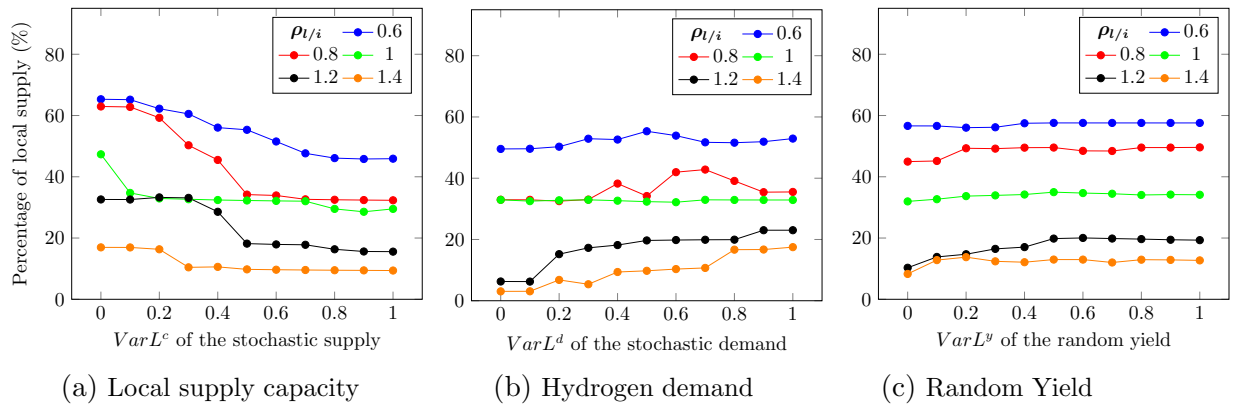
\begin{figure}[h]
\begin{subfigure}{0.33\textwidth}
    \centering
        \begin{tikzpicture}[scale=0.68, transform shape]
            \begin{axis}[
        xlabel={$VarL^c$ of the stochastic supply},
        ylabel={Percentage of local supply (\%) },
        ymin=0, ymax=95,
        xticklabel style={/pgf/number format/1000 sep=},
        legend pos=north east,
        ymajorgrids=false,
        grid style=dashed,
        legend columns=2,
    ]
    \addlegendimage{empty legend}
\addplot[blue,mark=*] coordinates {
    (0.0, 65.30589635)
    (0.1, 65.17150396)
    (0.2, 62.22750173)
    (0.3, 60.51982289)
    (0.4, 56.05624937)
    (0.5, 55.33330589)
    (0.6, 51.52074546)
    (0.7, 47.65170306)
    (0.8, 46.09077174)
    (0.9, 45.80675283)
    (1.0, 45.90275787)
};

\addplot[red,mark=*] coordinates {
    (0.0, 62.94733224)
    (0.1, 62.79353214)
    (0.2, 59.25443455)
    (0.3, 50.28713972)
    (0.4, 45.49522673)
    (0.5, 34.18256661)
    (0.6, 33.89951666)
    (0.7, 32.63692225)
    (0.8, 32.49432583)
    (0.9, 32.39550994)
    (1.0, 32.31931623)
};

\addplot[green,mark=*] coordinates {
    (0.0, 47.33684115)
    (0.1, 34.72662092)
    (0.2, 32.92065581)
    (0.3, 32.66273483)
    (0.4, 32.41434629)
    (0.5, 32.25545639)
    (0.6, 32.12292273)
    (0.7, 32.02712522)
    (0.8, 29.51342756)
    (0.9, 28.56906192)
    (1.0, 29.51342756)
};

\addplot[black,mark=*] coordinates {
    (0.0, 32.59356138)
    (0.1, 32.57032729)
    (0.2, 33.2660462)
    (0.3, 33.11163076)
    (0.4, 28.56906192)
    (0.5, 18.17528256)
    (0.6, 17.91514605)
    (0.7, 17.79450336)
    (0.8, 16.32746106)
    (0.9, 15.60724488)
    (1.0, 15.55740288)
};

\addplot[orange,mark=*] coordinates {
    (0.0, 16.95050004)
    (0.1, 16.92849808)
    (0.2, 16.32338787)
    (0.3, 10.43433693)
    (0.4, 10.58373487)
    (0.5, 9.770592689)
    (0.6, 9.655195044)
    (0.7, 9.563448734)
    (0.8, 9.487806654)
    (0.9, 9.438424996)
    (1.0, 9.389656606)
};

     \addlegendentry{\hspace{-.75cm}$\boldsymbol{\rho_{l/i}}$}
     \addlegendentry{0.6} 
      \addlegendentry{0.8} 
       \addlegendentry{1} 
        \addlegendentry{1.2} 
        \addlegendentry{1.4} 
          \end{axis}
        \end{tikzpicture}
    \caption{\small Local supply capacity}\label{fig:DifferentStocCap}
\end{subfigure}
\begin{subfigure}{0.32\textwidth}
    \centering
        \begin{tikzpicture}[scale=0.68, transform shape]
            \begin{axis}[
        xlabel={$VarL^d$ of the stochastic demand},
        ymin=0, ymax=95,
        xticklabel style={/pgf/number format/1000 sep=},
        legend pos=north east,
        ymajorgrids=false,
        grid style=dashed,
        legend columns=2,
    ]
    \addlegendimage{empty legend}
\addplot[blue,mark=*] coordinates {
    (0.0, 49.55898028)
    (0.1, 49.61458028)
    (0.2, 50.29011584)
    (0.3, 52.89476076)
    (0.4, 52.62181619)
    (0.5, 55.33330589)
    (0.6, 53.9023139)
    (0.7, 51.71705943)
    (0.8, 51.59150795)
    (0.9, 51.92435286)
    (1.0, 52.94025426)
};

\addplot[red,mark=*] coordinates {
    (0.0, 32.97557183)
    (0.1, 32.94501044)
    (0.2, 32.52068241)
    (0.3, 33.02556808)
    (0.4, 38.24932003)
    (0.5, 34.18256661)
    (0.6, 41.97357077)
    (0.7, 42.82111373)
    (0.8, 39.11610042)
    (0.9, 35.44496837)
    (1.0, 35.52066233)
};

\addplot[green,mark=*] coordinates {
    (0.0, 32.94501044)
    (0.1, 32.52068241)
    (0.2, 32.79860938)
    (0.3, 32.94779339)
    (0.4, 32.66273483)
    (0.5, 32.35423719)
    (0.6, 32.16892714)
    (0.7, 32.94401644)
    (0.8, 32.90778592)
    (0.9, 32.89207534)
    (1.0, 32.89207534)
};

\addplot[black,mark=*] coordinates {
    (0.0, 6.285916216)
    (0.1, 6.254036462)
    (0.2, 15.20906651)
    (0.3, 17.26824518)
    (0.4, 18.17528256)
    (0.5, 19.68405011)
    (0.6, 19.81801641)
    (0.7, 19.90592085)
    (0.8, 19.93812804)
    (0.9, 23.02324179)
    (1.0, 23.02324179)
};

\addplot[orange,mark=*] coordinates {
    (0.0, 3.060676247)
    (0.1, 3.081465556)
    (0.2, 6.797312115)
    (0.3, 5.381718677)
    (0.4, 9.38027969)
    (0.5, 9.770592689)
    (0.6, 10.33890192)
    (0.7, 10.71112057)
    (0.8, 16.67294906)
    (0.9, 16.71362565)
    (1.0, 17.49799295)
};

     \addlegendentry{\hspace{-.75cm}$\boldsymbol{\rho_{l/i}}$}
     \addlegendentry{0.6} 
      \addlegendentry{0.8} 
       \addlegendentry{1} 
        \addlegendentry{1.2} 
        \addlegendentry{1.4} 
          \end{axis}
        \end{tikzpicture}
    \caption{\small Hydrogen demand}\label{fig:DifferentDemand}
\end{subfigure}
\begin{subfigure}{0.25\textwidth}
    \centering
        \begin{tikzpicture}[scale=0.68, transform shape]
            \begin{axis}[
        xlabel={$VarL^y$ of the random yield},
        ymin=0, ymax=95,
        xticklabel style={/pgf/number format/1000 sep=},
        legend pos=north east,
        ymajorgrids=false,
        grid style=dashed,
        legend columns=2,
    ]
    \addlegendimage{empty legend}
\addplot[blue,mark=*] coordinates {
    (0.0, 56.61189458)
    (0.1, 56.59809055)
    (0.2, 56.03385789)
    (0.3, 56.15896651)
    (0.4, 57.46837194)
    (0.5, 57.58250275)
    (0.6, 57.57911105)
    (0.7, 57.57887116)
    (0.8, 57.57906147)
    (0.9, 57.57399121)
    (1.0, 57.57906147)
};

\addplot[red,mark=*] coordinates {
    (0.0, 44.95088421)
    (0.1, 45.15973148)
    (0.2, 49.30366242)
    (0.3, 49.18412376)
    (0.4, 49.50481441)
    (0.5, 49.53807252)
    (0.6, 48.45958166)
    (0.7, 48.41435594)
    (0.8, 49.53485252)
    (0.9, 49.53407252)
    (1.0, 49.59804852)
};

\addplot[green,mark=*] coordinates {
    (0.0, 31.95632486)
    (0.1, 32.67631918)
    (0.2, 33.68476982)
    (0.3, 33.92213449)
    (0.4, 34.20280885)
    (0.5, 34.97874573)
    (0.6, 34.68385522)
    (0.7, 34.46216479)
    (0.8, 34.05613199)
    (0.9, 34.18852875)
    (1.0, 34.12331035)
};

\addplot[black,mark=*] coordinates {
    (0.0, 10.29452208)
    (0.1, 13.82571471)
    (0.2, 14.70898364)
    (0.3, 16.45643892)
    (0.4, 17.03609635)
    (0.5, 19.81932662)
    (0.6, 20.03148849)
    (0.7, 19.83432665)
    (0.8, 19.64811338)
    (0.9, 19.43594561)
    (1.0, 19.30911012)
};

\addplot[orange,mark=*] coordinates {
    (0.0, 8.276600627)
    (0.1, 12.79628195)
    (0.2, 13.76103976)
    (0.3, 12.39252954)
    (0.4, 12.10599849)
    (0.5, 12.97104594)
    (0.6, 12.959450837)
    (0.7, 12.01522189)
    (0.8, 12.92935926)
    (0.9, 12.86177255)
    (1.0, 12.70304946)
};
     \addlegendentry{\hspace{-.75cm}$\boldsymbol{\rho_{l/i}}$}
     \addlegendentry{0.6} 
      \addlegendentry{0.8} 
       \addlegendentry{1} 
        \addlegendentry{1.2} 
        \addlegendentry{1.4} 
          \end{axis}
        \end{tikzpicture}
    \caption{\small Random Yield}\label{fig:DifferentRandomYield}
\end{subfigure}

\caption{Impact of the variability levels on the local supply rate for various cost ratios}\label{fig:DifferentVariability}
\end{figure}

In this section, we analyze the variability level of local supply capacity ($VarL^c$), demand ($VarL^d$), and random yield ($VarL^y$), holding the mean of their distributions constant while adjusting the standard deviation.  Understanding the impact of these uncertainties is crucial to determine the extent to which each source of variability influences the reliance on local production versus import. Figure~\ref{fig:DifferentVariability} provides the results under different variability levels and cost ratios ($\rho_{l/i}$). 

In Figure~\ref{fig:DifferentStocCap}, we keep the variability of the demand and the random yield constant at 0.5 and analyze the variability of the local supply capacity. The graph shows that as uncertainty in supply increases, the system relies less on local production and more on import. This is expected as a less reliable local supply forces the system to shift towards importing to mitigate risk. We observe that changes in supply variability have a significant impact on the percentage of local supply. For example, when $\rho_{l/i}$ is 0.8, increasing the level of variability results in a reduction of up to 30\% in local supply. These results underscore the critical need to account for stochastic supply in strategic planning. In real-world scenarios, if local renewable energy sources, such as wind and solar power, become highly unpredictable because of adverse weather conditions, decision-makers may need to rely more on importing to maintain stable hydrogen supply and meet demand without significant disruptions.
 
In Figure~\ref{fig:DifferentDemand}, we fix the local supply capacity and the random yield variability to 0.5 and change the demand variability. We observe that as the variability of demand increases, the system increasingly relies on local production. This is due to local supplier’s ability to respond more quickly to fluctuations compared to importing, which have longer lead times. For example, in real-world scenarios, if demand in the hydrogen market becomes more unpredictable due to dynamic industrial activities, shifting market needs, or sudden changes in policy incentives, decision-makers can rely more on local production. These significant changes in the percentage of local supply show that accounting for demand uncertainty is also critical in the decision-making process.

In Figure~\ref{fig:DifferentRandomYield}, we keep the supply capacity and demand variability constant at 0.5 and analyze the impact of variability on the random yield. We observe that this variability does not significantly affect the percentage of local supply when local supply is less costly than importing. However, when importing is cheaper than local production, such as at $\rho_{l/i}$ values of 1.2 and 1.4, increasing random yield variability leads to a shift towards local supply. In practical terms, when the portion of hydrogen lost during transportation and conversion processes becomes highly unpredictable due to logistical inefficiencies, weather-related disruptions, or technological limitations, the sourcing policy often shifts to relying more on local production to maintain a stable and reliable supply. This highlights the importance of considering random yield in the model, particularly when importing is a key part of the policy. In this analysis, we set the average random yield at 17.5\%, but as this average increases, its influence on the system could become even more pronounced.

\section{Conclusions}\label{sec:concluions}
In this paper, we study a dual sourcing inventory control model for green hydrogen, accounting for general lead times, stochastic local supply capacity, and random yield from import. We formulate the problem as a Markov decision process and solve it to optimality using value iteration. In addition to the optimal dual sourcing policy, we propose heuristic policies designed to stabilize order quantities, which is appealing to practitioners in forming hydrogen trade agreements within complex supply chain dynamics.

We build a case study for the Netherlands. We show that our dual sourcing model, which considers stochastic supply capacity and random yield, outperforms its counterparts, such as single-sourcing models and dual sourcing models that ignore stochastic capacity, random yield, or both. We show that considering both stochastic supply capacity and random yield demonstrates an average cost benefit of 8\%. We evaluate the performance of the proposed heuristic policies \texttt{FOQ}, \texttt{FOQ+}, \texttt{TBS}, and \texttt{TBS+} in comparison to the optimal policy, identifying conditions under which they achieve comparable results. For example, our experiments show that transitioning from the stable \texttt{FOQ} policy to the partially adjustable \texttt{FOQ+} policy improves cost efficiency by an average of 11\% and performs within 2\% of the optimal policy. In turn, our findings can help decision-makers to shape hydrogen trade agreements between importing and exporting countries. Furthermore, our sensitivity analyses identify the necessary conditions—such as production and storage costs, as well as variability in local supply capacity, hydrogen demand, and random yield—to achieve specific local supply levels. These analyses support policymakers in evaluating the feasibility of the Netherlands' climate scenarios, which target different levels of local supply and imports, and in guiding strategic investment planning.

Future research can focus on developing tailored solution approaches, such as deep reinforcement learning algorithms, to address large-scale systems. Nevertheless, a key question concerns how such advanced policies can ensure order stability in complex supply chain environments. Achieving stability with these approaches remains an open research question that requires further investigation.
In addition, it would be interesting to study the integration of national-level hydrogen dual sourcing with national distribution networks, either based on truck or pipeline transportation. Special emphasis can be placed on geographically strategic storage and demand locations and on the logistics of transporting hydrogen between these points.

\section{Acknowledgements}
This project has received funding from the Fuel Cells and Hydrogen 2 Joint Undertaking (now Clean Hydrogen Partnership) under Grant Agreement No 875090. This Joint Undertaking receives support from the European Union's Horizon 2020 research and innovation programme, Hydrogen Europe and Hydrogen Europe Research. Albert H. Schrotenboer is supported by a VENI research talent grant from the Dutch Science Foundation (NWO).

\bibliographystyle{apalike}
\bibliography{references}
\end{document}